\documentclass[12pt]{article}
\usepackage{amssymb}

\newtheorem{definition}{Definition}[section]
\newtheorem{theorem}[definition]{Theorem}
\newtheorem{lemma}[definition]{Lemma}

\newtheorem{remark}[definition]{Remark}
\newtheorem{example}[definition]{Example}

\newtheorem{problem}[definition]{Problem}
\newtheorem{note}[definition]{Note}

\def\C{\mathbb C}
\def\K{\mathbb K}

\def\K{\mathbb K}

\begin{document}

\title{ \bf Tridiagonal pairs and the quantum\\
affine algebra $U_q({\widehat {sl}}_2)$}
\author{Tatsuro Ito and Paul Terwilliger
}
\date{}
\maketitle
\begin{abstract} 
Let $\K$ denote an algebraically closed field and
let $q$ denote a nonzero scalar in $\K$ that is not a root of
unity.
Let $V$ denote a vector space over $\K$ with finite positive
dimension and
let $A,A^*$ denote a tridiagonal pair on $V$.
Let $\theta_0, \theta_1,\ldots, \theta_d$
(resp. 
$\theta^*_0, \theta^*_1,\ldots, \theta^*_d$)
denote a standard ordering of the eigenvalues of $A$
(resp.  $A^*$).
We assume 
there exist nonzero
scalars $a, a^*$ in $\K$ 
such that
$\theta_i =  aq^{2i-d}$
and
$\theta^*_i =  a^*q^{d-2i}$
for $0 \leq i \leq d$.
We display two irreducible 
$U_q({\widehat{sl}}_2)$-module structures on $V$ 
and discuss how these are related to the actions of $A$ and $A^*$. 
\end{abstract}

\section{The quantum affine algebra 
$U_q(\widehat{ sl}_2)$ }

\medskip
\noindent Throughout this paper  
 $\K$ will  denote an algebraically closed field. We fix
 a nonzero scalar
 $q\in \K$ that is not
 a root of unity.
We will use the following notation.

\begin{equation}
\lbrack n \rbrack_q = {{q^n-q^{-n}}\over {q-q^{-1}}},
\qquad \qquad n=0,1,\ldots
\label{eq:brackndef}
\end{equation}

\noindent We now recall the definition of
$U_q(\widehat{ sl}_2)$.

\begin{definition} 
\label{def:uq}
\rm
\cite[p.~262]{cp}
The quantum affine algebra 
$U_q(\widehat{ sl}_2)$ is the unital associative $\K$-algebra 
with
generators $e^{\pm}_i$, $K_i^{{\pm}1}$, $i\in \lbrace 0,1\rbrace $
and the following relations:
\begin{eqnarray}
K_iK^{-1}_i &=& 
K^{-1}_iK_i =  1,
\label{eq:buq1}
\\
K_0K_1&=& K_1K_0,
\label{eq:buq2}
\\
K_ie^{\pm}_iK^{-1}_i &=& q^{{\pm}2}e^{\pm}_i,
\label{eq:buq3}
\\
K_ie^{\pm}_jK^{-1}_i &=& q^{{\mp}2}e^{\pm}_j, \qquad i\not=j,
\label{eq:buq4}
\\
\lbrack e^+_i, e^-_i\rbrack &=& {{K_i-K^{-1}_i}\over {q-q^{-1}}},
\label{eq:buq5}
\\
\lbrack e^{\pm}_0, e^{\mp}_1\rbrack &=& 0,
\label{eq:buq6}
\end{eqnarray}
\begin{eqnarray}
(e^{\pm}_i)^3e^{\pm}_j -  
\lbrack 3 \rbrack_q (e^{\pm}_i)^2e^{\pm}_j e^{\pm}_i 
+\lbrack 3 \rbrack_q e^{\pm}_ie^{\pm}_j (e^{\pm}_i)^2 - 
e^{\pm}_j (e^{\pm}_i)^3 =0, \qquad i\not=j.
\label{eq:buq7}
\end{eqnarray}
We call $e^{\pm}_i$, $K_i^{{\pm}1}$, $i\in \lbrace 0,1\rbrace $
the {\it Chevalley generators} for
$U_q({\widehat{sl}}_2)$.

\end{definition}

\begin{remark}
\rm The equations
(\ref{eq:buq7})
are called the
{\it $q$-Serre relations}.
\end{remark}

\section{A presentation of 
$U_q(\widehat{ sl}_2)$}

In order to state our main result 
we introduce an alternate presentation of
$U_q(\widehat{ sl}_2)$. This presentation is given below.

\begin{theorem}
\label{thm:uq2}
The quantum affine algebra
$U_q(\widehat{ sl}_2)$ is isomorphic to
the unital associative $\K$-algebra 
with
generators $y^{\pm}_i$, $k_i^{{\pm}1}$, $i\in \lbrace 0,1\rbrace $
and the following relations:
\begin{eqnarray}
k_ik^{-1}_i = 
k^{-1}_ik_i &=&  1,
\label{eq:2buq1}
\\
k_0k_1 \; \mbox{is central},
\label{eq:2buq2}
\\
\frac{q y^+_ik_i-q^{-1}k_iy^+_i}{q-q^{-1}} &=& 1,
\label{eq:2buq3}
\\
\frac{q k_iy^-_i-q^{-1}y^-_ik_i}{q-q^{-1}} &=& 1,
\label{eq:2buq4}
\\
\frac{q y^-_iy^+_i-q^{-1}y^+_iy^-_i}{q-q^{-1}} &=& 1,
\label{eq:2buq5}
\\
\frac{q y^+_iy^-_j-q^{-1}y^-_jy^+_i}{q-q^{-1}} &=& k^{-1}_0k^{-1}_1,
\qquad i\not=j,
\label{eq:2buq6}
\end{eqnarray}
\begin{eqnarray}
(y^{\pm}_i)^3y^{\pm}_j -  
\lbrack 3 \rbrack_q (y^{\pm}_i)^2y^{\pm}_j y^{\pm}_i 
+\lbrack 3 \rbrack_q y^{\pm}_iy^{\pm}_j (y^{\pm}_i)^2 - 
y^{\pm}_j (y^{\pm}_i)^3 =0, \qquad i\not=j.
\label{eq:2buq7}
\end{eqnarray}
An isomorphism with the presentation in Definition
\ref{def:uq} is given by:
\begin{eqnarray*}
\label{eq:iso1}
k^\pm_i &\rightarrow & K^\pm_i,\\
\label{eq:iso2}
y^-_i &\rightarrow & K^{-1}_i+e^{-}_i, \\
\label{eq:iso3}
y^+_i &\rightarrow & K^{-1}_i-q(q-q^{-1})^2K^{-1}_ie^{+}_i.
\end{eqnarray*}
The inverse of this isomorphism is given by:
\begin{eqnarray*}
\label{eq:iso1inv}
K^\pm_i &\rightarrow & k^\pm_i,\\
\label{eq:iso2inv}
e^-_i &\rightarrow & y^{-}_i-k^{-1}_i, \\
\label{eq:iso3inv}
e^+_i &\rightarrow & \frac{1-k_i y^{+}_i}{q(q-q^{-1})^2}.
\end{eqnarray*}
\end{theorem}
\noindent {\it Proof:} One readily checks that
each map is a homomorphism of $\K$-algebras and that
the maps are inverses.
It follows each map is an isomorphism of $\K$-algebras.
\hfill $\Box $ \\

\begin{definition}
\label{def:alt}
\rm 
With reference to Theorem
\ref{thm:uq2}
we call
 $y^{\pm}_i$, $k_i^{{\pm}1}$, $i\in \lbrace 0,1\rbrace $
the {\it alternate generators} of
$U_q({\widehat{sl}}_2)$.
\end{definition}

\section{Tridiagonal pairs}
\noindent 
We now recall the notion of a {\it tridiagonal pair}
\cite{TD00}, \cite{qSerre}.
We will use the following terms.
 Let $V$  denote
a vector space over $\K$ with finite positive dimension.
 Let $A:V\rightarrow V$ denote a linear transformation and
 let $W$ denote a subspace of $V$. We call $W$ an {\it eigenspace} of $A$ whenever $W\not=0$ and there exists $\theta \in \K$ such that 
\begin{eqnarray*}
W=\lbrace v \in V \;\vert \;Av = \theta v\rbrace.
\end{eqnarray*}
We say $A$ is {\it diagonalizable} whenever
$V$ is spanned by the eigenspaces of $A$.

\begin{definition}  
\cite[Definition 1.1]{TD00}
\label{def:tdp}
\rm
Let $V$ denote a vector space over $\K$ with finite positive dimension.
By a {\it tridiagonal pair}  on $V$,
we mean an ordered pair $A,A^*$ where
$A:V\rightarrow V$ and 
$A^*:V\rightarrow V$ 
 are linear transformations that satisfy
the following four conditions.
\begin{enumerate}
\item Each of $A,A^*$ is diagonalizable.
\item There exists an ordering $V_0, V_1,\ldots, V_d$ of the  
eigenspaces of $A$ such that 
\begin{equation}
A^* V_i \subseteq V_{i-1} + V_i+ V_{i+1} \qquad \qquad (0 \leq i \leq d),
\label{eq:t1}
\end{equation}
where $V_{-1} = 0$, $V_{d+1}= 0$.
\item There exists an ordering $V^*_0, V^*_1,\ldots, V^*_\delta$ of
the  
eigenspaces of $A^*$ such that 
\begin{equation}
A V^*_i \subseteq V^*_{i-1} + V^*_i+ V^*_{i+1} \qquad \qquad (0 \leq i \leq \delta),
\label{eq:t2}
\end{equation}
where $V^*_{-1} = 0$, $V^*_{\delta+1}= 0$.
\item There does not exist a subspace $W$ of $V$ such  that $AW\subseteq W$,
$A^*W\subseteq W$, $W\not=0$, $W\not=V$.
\end{enumerate}
\end{definition}

\begin{note} \rm
According to a common 
notational convention, $A^*$ denotes the conjugate transpose
of $A$.
 We are not using this convention.
In a tridiagonal pair $A,A^*$ the linear transformations 
$A$ and $A^*$ are arbitrary subject to (i)--(iv) above.
\end{note}

\noindent 
Our interest in tridiagonal pairs evolved from our interest in
the 
following special case.
A tridiagonal pair for which the $V_i, V^*_i$ all
have dimension 1 is called a {\it Leonard pair}
\cite{LS99}.
There is a natural correspondence between
the Leonard pairs and a family of orthogonal
polynomials consisting of the $q$-Racah polynomials
\cite{AWil},
\cite{GR}
and some related polynomials in the Askey-scheme
\cite{KoeSwa},
\cite{TLT:array}.
This correspondence follows from
the classification of Leonard pairs
\cite{LS99},
\cite{TLT:array}.
We remark that 
this classification amounts to a linear algebraic
version of a theorem of D. Leonard
\cite{BanIto},
\cite{Leon}
concerning the $q$-Racah polynomials.
See 
\cite{shape},
\cite{qSerre},
   \cite{LS24},
   \cite{conform},
    \cite{lsint},
\cite{Terint},
 \cite{TLT:split},
\cite{qrac}, 
\cite{aw}
for more information about Leonard pairs.

\medskip
\noindent 
Given these comments on Leonard pairs,
it is natural to attempt a classification of
the 
tridiagonal pairs.
At present we do not have this classification;
however we do have 
a result that might lead to
one.
In order to state the result we
recall a few basic facts about tridiagonal pairs.
Let $A,A^*$ denote a tridiagonal pair on $V$
and let $d, \delta$ be as in
Definition \ref{def:tdp}(ii), (iii).
 By 
 \cite[Lemma 4.5]{TD00}
 we have
$d=\delta$; we call this common value the {\it diameter}
of $A,A^*$.  
An ordering of the eigenspaces of $A$ (resp. $A^*$)
will be called
{\it standard} whenever it satisfies
(\ref{eq:t1}) (resp. 
(\ref{eq:t2})). 
We comment on the uniqueness of the standard ordering.
Let 
 $V_0, V_1,\ldots, V_d$ denote a standard ordering
of the eigenspaces of $A$.
Then the ordering 
 $V_d, V_{d-1},\ldots, V_0$ 
is standard and no other ordering is standard.
A similar result holds for the 
 eigenspaces of $A^*$.
An ordering of the eigenvalues of $A$ (resp. $A^*$)
will be called {\it standard}
whenever the corresponding ordering of the eigenspaces of $A$
(resp. $A^*$) is standard.
Let $\theta_0,\theta_1,\ldots, \theta_d$ (resp. 
 $\theta^*_0,\theta^*_1,\ldots, \theta^*_d$)
 denote a standard ordering
of the eigenvalues of $A$ (resp. $A^*$).
The $\theta_i, \theta^*_i$ satisfy a number of
equations 
\cite[Theorem 4.3]{qSerre} that have been solved in closed
form \cite[Theorem 4.4]{qSerre}.
In a special case of interest,
there exist nonzero scalars $a, a^*$ in $\K$ such
that $\theta_i = aq^{2i-d}$ and $\theta^*_i=a^*q^{d-2i}$
for $0 \leq i \leq d$ 
 \cite[Example 1.7]{TD00}, \cite{shape}.

\medskip
\noindent We now state our main result.

\begin{theorem}
\label{thm:main}
Let $V$ denote a vector space over $\K$ with finite positive
dimension and
let $A,A^*$ denote a tridiagonal pair on $V$.
Let $\theta_0, \theta_1,\ldots, \theta_d$
(resp. 
$\theta^*_0, \theta^*_1,\ldots, \theta^*_d$)
denote a standard ordering of the eigenvalues of $A$
(resp.  $A^*$).
We assume 
there exist nonzero
scalars $a, a^*$ in $\K$ such that
$\theta_i =  aq^{2i-d}$
and
$\theta^*_i =  a^*q^{d-2i}$
for $0 \leq i \leq d$.
Then with reference to Theorem  
\ref{thm:uq2},
there exists a unique 
$U_q({\widehat{sl}}_2)$-module structure on $V$ 
such that $ay_1^-$ acts as $A$ and $a^*y_0^-$
acts as $A^*$. 
Moreover there exists
a unique $U_q({\widehat{sl}}_2)$-module structure on $V$ 
such that $a y_0^+$ acts as $A$ and $a^*y_1^+$
acts as $A^*$. 
Both $U_q({\widehat{sl}}_2)$-module
 structures are 
irreducible.
\end{theorem}

\noindent The proof of Theorem
\ref{thm:main} appears in Sections 13, 14 below.

\begin{remark}
\rm
The finite dimensional irreducible 
modules for
$U_q({\widehat{sl}}_2)$ are described in
\cite{cp}. In a future paper we hope 
to use 
\cite{cp}
to obtain a classification of the
tridiagonal pairs that satisfy the assumptions of
Theorem
\ref{thm:main}. See Lemma 
\ref{lem:whenir} and Problem
\ref{prob:one} below for a discussion of the issues involved.
\end{remark}

\begin{remark}
\rm 
Theorem \ref{thm:main} extends some work
of Curtin and Al-Najjar \cite{ch}, \cite{ch2}.
They
give
a
 $U_q({\widehat{sl}}_2)$-action
for those
tridiagonal pairs that satisfy the
assumptions of 
Theorem \ref{thm:main} 
and for which the dimensions of the $V_i, V^*_i$ are
all at most 2.
\end{remark}

\section{Six decompositions}
In this section and the next we collect some results about
tridiagonal pairs which we will use to prove
Theorem
\ref{thm:main}.

\medskip
\noindent We will use the following notation.
Let $V$ denote a vector space over $\K$ with  finite positive
dimension.
Let $d$ denote a nonnegative integer.
By a {\it decomposition of $V$ of length $d$}, we mean
a sequence $U_0, U_1, \ldots, U_d$ consisting of nonzero
subspaces of $V$ such that
\begin{eqnarray*}
\label{eq:dec}
V = U_0+ U_1 + \cdots + U_d \qquad (\mbox{direct sum}).
\end{eqnarray*}
We do not assume each of $U_0, U_1, \ldots, U_d$ has  dimension 1.
For $0 \leq i \leq d$ we call $U_i$ the {\it $i$th subspace} of 
the decomposition.
For notational convenience we define $U_{-1}:=0$ and $U_{d+1}:=0$.

\medskip
\noindent
We will refer to the following setup.

\begin{definition}
\label{def:setup}
\rm
Let $V$ denote a vector space over $\K$ with finite positive
dimension  and let $A,A^*$ denote a tridiagonal pair on $V$.
Let $V_0, V_1, \ldots, V_d$ (resp.  
$V^*_0, V^*_1, \ldots, V^*_d$)
denote a standard ordering of the eigenspaces of $A$
(resp. $A^*$).
For $0 \leq i\leq d$ let $\theta_i$ (resp. $\theta^*_i$)
denote the eigenvalue of $A$ (resp. $A^*$) associated
with $V_i$ (resp. $V^*_i$).
\end{definition}

\noindent With reference to Definition
\ref{def:setup},
we are about to define
six decompositions of $V$.
In order to keep track of these decompositions
 we will give each of them
a name. Our naming scheme is as follows.
Let $\Omega$ denote the set
consisting of the four symbols
$0, D, 0^*, D^*$. Each of the six decompositions will
get a name $\lbrack u \rbrack$ where $u$ is a two-element
subset of $\Omega$.
We now define the six decompositions.

\begin{lemma}
\label{thm:sixdecp}
With reference to Definition
\ref{def:setup},
for each of the six rows in the table below,
and for $0 \leq i \leq d$, let $U_i$ denote
the ith subspace  described in that row.
Then the sequence $U_0, U_1, \ldots, U_d$ is a decomposition
of $V$.
\medskip

\centerline{
\begin{tabular}[t]{c|c}
       {\rm name} & {\rm $i$th subspace of the decomposition }
 \\ \hline  \hline
	$\lbrack 0D\rbrack$ & $V_i$    
	\\
	$\lbrack 0^*D^*\rbrack$ & $V^*_i$   \\
	$\lbrack 0^*D\rbrack$ &
	$(V^*_0+\cdots + V^*_i)\cap (V_i+\cdots +V_d)$   \\ 
	$\lbrack 0^*0\rbrack$ & 
	$(V^*_{0}+\cdots +V^*_i) \cap
	(V_0+\cdots +V_{d-i})$
	\\ 
        $\lbrack D^*0\rbrack$ &
	$
	(V^*_{d-i}+\cdots +V^*_{d}) \cap
	(V_0+\cdots + V_{d-i})$
	\\ 
	$ \lbrack D^*D\rbrack $ & 
	$
	(V^*_{d-i}+\cdots +V^*_{d})
	\cap
	(V_i+\cdots +V_{d})
	$ 
	\end{tabular}}
\medskip
\noindent 
\end{lemma}
\noindent {\it Proof:} We consider each of the six rows
of the table.
\\
\noindent $\lbrack 0D\rbrack$: Recall $V_0, V_1, \ldots, V_d$
are the eigenspaces of $A$ and that $A$ is diagonalizable.
\\
\noindent $\lbrack 0^*D^*\rbrack$: Recall $V^*_0, V^*_1, \ldots, V^*_d$
are the eigenspaces of $A^*$ and that $A^*$ is diagonalizable.
\\
\noindent $\lbrack 0^*D\rbrack$: 
Define $U_i=
	(V^*_0+\cdots + V^*_i)\cap (V_i+\cdots +V_d)$ 
for $0 \leq i \leq d$. Then the sequence $U_0, U_1, \ldots, U_d$
is a decomposition of $V$ by \cite[Theorem 4.6]{TD00}.
\\
\noindent $\lbrack 0^*0\rbrack$: 
Apply the present Lemma, row
 $\lbrack 0^*D\rbrack$, with $V_i$ replaced by $V_{d-i}$ for
$0 \leq i \leq d$.
\\
\noindent $\lbrack D^*0\rbrack$: 
Apply the present Lemma, row
 $\lbrack 0^*D\rbrack$, with
 $V_i$ replaced by $V_{d-i}$ and
 $V^*_i$ replaced by $V^*_{d-i}$
 for
$0 \leq i \leq d$.
\\
\noindent $\lbrack D^*D\rbrack$: 
Apply the present Lemma, row
 $\lbrack 0^*D\rbrack$, with
 $V^*_i$ replaced by $V^*_{d-i}$
 for
$0 \leq i \leq d$.
\hfill $\Box $ \\

\noindent The six decompositions from
Lemma
\ref{thm:sixdecp} are related to each other as follows.

\begin{lemma}
\label{thm:decsum}
Adopt the assumptions of Definition
\ref{def:setup} and
let $U_0, U_1, \ldots, U_d$ 
denote any one of the six decompositions
of $V$ given in 
Lemma \ref{thm:sixdecp}.
Then for $0 \leq i \leq d$
the sums $U_0+\cdots + U_i$ and $U_i+\cdots + U_d$
are given as follows.

\medskip

\centerline{
\begin{tabular}[t]{c|c|c}
      {\rm name} &$U_0+\cdots + U_i$ & $U_i+\cdots + U_d$ \\ \hline  \hline
	$\lbrack 0D\rbrack $ &
	$V_0 + \cdots + V_i$  & $V_i + \cdots + V_d $  
\\
$\lbrack 0^*D^*\rbrack $ &
	$V^*_0+\cdots + V^*_i$  & $ V^*_i + \cdots + V^*_d$   \\
       $\lbrack 0^*D\rbrack $ &
        $V^*_0+\cdots + V^*_i$ & $V_i+\cdots +V_d$   \\ 
	$\lbrack 0^*0\rbrack $ & 
	$V^*_{0}+\cdots +V^*_i$  
	&
	$V_{0}+\cdots +V_{d-i}$
	\\ 
       $\lbrack D^*0\rbrack $ & 
	$V^*_{d-i}+\cdots + V^*_{d}$ 
&	
	$V_{0} + \cdots +V_{d-i} $ 
	\\ 
	$\lbrack D^*D\rbrack $ &
	$V^*_{d-i}+\cdots +V^*_{d}$ 
&	
	$V_i+\cdots +V_{d} $
	\end{tabular}}
\medskip

\end{lemma}
\noindent {\it Proof:} We consider each of the six rows of the table.
\\
\noindent $\lbrack 0D\rbrack$:
Immediate from Lemma
 \ref{thm:sixdecp}, row
 $\lbrack 0D\rbrack$.
\\
\noindent $\lbrack 0^*D^*\rbrack$: 
Immediate from Lemma
 \ref{thm:sixdecp}, row
 $\lbrack 0^*D^*\rbrack$.
\\
\noindent $\lbrack 0^*D\rbrack$: 
Let $U_0, U_1, \ldots, U_d$ denote the decomposition
$\lbrack 0^*D\rbrack$. By \cite[Theorem 4.6]{TD00}
we find $U_0+\cdots + U_i=V^*_0+\cdots + V^*_i$
and 
$U_i+\cdots + U_d=V_i+\cdots + V_d$ for $0 \leq i \leq d$.
\\
\noindent $\lbrack 0^*0\rbrack$: 
Apply the present Lemma, row
 $\lbrack 0^*D\rbrack$, with $V_i$ replaced by $V_{d-i}$ for
$0 \leq i \leq d$.
\\
\noindent $\lbrack D^*0\rbrack$: 
Apply the present Lemma, row
 $\lbrack 0^*D\rbrack$, with
 $V_i$ replaced by $V_{d-i}$ and
 $V^*_i$ replaced by $V^*_{d-i}$
 for
$0 \leq i \leq d$.
\\
\noindent $\lbrack D^*D\rbrack$: 
Apply the present Lemma, row
 $\lbrack 0^*D\rbrack$, with
 $V^*_i$ replaced by $V^*_{d-i}$
 for
$0 \leq i \leq d$.
\hfill $\Box $ \\

\noindent We have a comment.

\begin{lemma}
\cite[Corollary 5.7, Corollary 6.6]{TD00}
\label{lem:shape}
Adopt the assumptions of Definition
\ref{def:setup} and 
let $U_0, U_1, $ $ \ldots, U_d$ 
denote any one of the six decompositions
of $V$ given in 
Lemma \ref{thm:sixdecp}.
For $0 \leq i \leq d$ let $\rho_i$ denote the dimenension
of $U_i$. Then the sequence
$\rho_0, \rho_1, \ldots, \rho_d$ is independent of
the decomposition.
Moreover the sequence $\rho_0, \rho_1, \ldots, \rho_d$
is unimodal and symmetric; that is
$ \rho_i =\rho_{d-i}$ for $0 \leq i \leq d$
and 
$\rho_{i-1} \leq \rho_{i}$ for  $1 \leq i \leq d/2$.
\end{lemma}

\noindent Referring to Lemma 
\ref{lem:shape}, we call the sequence $\rho_0, \rho_1, \ldots, \rho_d$
the {\it shape} of the tridiagonal pair.
As we indicated in Section 2,
a tridiagonal 
pair of shape $1,1,\ldots, 1$ is the same thing as a Leonard pair
\cite{TD00}.

\section{The action of $A$ and $A^*$ on the six decompositions}

\noindent With reference to Definition \ref{def:setup},
 in this section we describe the actions of
$A$ and $A^*$ on each of the six decompositions given
in Lemma \ref{thm:sixdecp}.

\begin{lemma}
\label{thm:aaction}
Adopt the assumptions of Definition \ref{def:setup} and
let $U_0, U_1, \ldots, U_d$ 
denote any one of the six decompositions
of $V$ given in 
Lemma \ref{thm:sixdecp}.
Then for $0 \leq i \leq d$
the action of $A$ and $A^*$ on $U_i$ is described as follows.

\medskip

\centerline{
\begin{tabular}[t]{c|c|c}
      {\rm name} &{\rm action of $A$ on $U_i$} & {\rm action of $A^*$ on $U_i$}
      \\ \hline  \hline
	$\lbrack 0D\rbrack $ &
	$ (A-\theta_iI)U_i=0$  & 
	$A^* U_i \subseteq U_{i-1}+ U_i+U_{i+1}$ 
\\	
	$\lbrack 0^*D^*\rbrack $ &
	$A U_i \subseteq U_{i-1}+ U_i+U_{i+1}$  &
	$ (A^*-\theta^*_iI)U_i=0$   \\
       $\lbrack 0^*D\rbrack $ &
        $(A-\theta_iI)U_i\subseteq U_{i+1}$ &
	$(A^*-\theta^*_iI)U_i \subseteq U_{i-1}$   \\ 
	$\lbrack 0^*0\rbrack $ & 
	$(A-\theta_{d-i}I)U_i\subseteq U_{i+1}$ &
	$(A^*-\theta^*_iI)U_i \subseteq U_{i-1}$   \\ 
       $\lbrack D^*0\rbrack $ & 
	$(A-\theta_{d-i}I)U_i\subseteq U_{i+1}$ &
	$(A^*-\theta^*_{d-i}I)U_i \subseteq U_{i-1}$   \\ 
	$\lbrack D^*D\rbrack $ &
	$(A-\theta_iI)U_i\subseteq U_{i+1}$ &
	$(A^*-\theta^*_{d-i}I)U_i \subseteq U_{i-1}$  
	\end{tabular}}
\medskip

\end{lemma}
\noindent {\it Proof:} We consider each of the six rows of the table.
\\
\noindent $\lbrack 0D\rbrack$: 
For $0\leq i \leq d$ the space
$V_i$ is an eigenspace for $A$ with
eigenvalue $\theta_i$. Therefore $(A-\theta_iI)V_i=0$.
We have 
$A^*V_i\subseteq V_{i-1}
+V_i
+ V_{i+1}$
by (\ref{eq:t1}).
\\
\noindent $\lbrack 0^*D^*\rbrack$:
For $0 \leq i \leq d$ 
we find
$AV^*_i\subseteq V^*_{i-1}
+V^*_i
+ V^*_{i+1}$
by (\ref{eq:t2}).
The space $V^*_i$ is an eigenspace for $A^*$ with
eigenvalue $\theta^*_i$. Therefore $(A^*-\theta^*_iI)V^*_i=0$.
\\
\noindent $\lbrack 0^*D\rbrack$: 
Let $U_0, U_1, \ldots, U_d$ denote the decomposition
 $\lbrack 0^*D\rbrack$.
By
 \cite[Theorem 4.6]{TD00} we find
$(A-\theta_iI)U_i\subseteq U_{i+1}$ and
$(A^*-\theta^*_iI)U_i\subseteq U_{i-1}$
for $0 \leq i \leq d$.
\\
\noindent $\lbrack 0^*0\rbrack$: 
Apply the present Lemma, row
 $\lbrack 0^*D\rbrack$, with $V_i$ replaced by $V_{d-i}$ for
$0 \leq i \leq d$.
\\
\noindent $\lbrack D^*0\rbrack$: 
Apply the present Lemma, row
 $\lbrack 0^*D\rbrack$, with
 $V_i$ replaced by $V_{d-i}$ and
 $V^*_i$ replaced by $V^*_{d-i}$
 for
$0 \leq i \leq d$.
\\
\noindent $\lbrack D^*D\rbrack$: 
Apply the present Lemma, row
 $\lbrack 0^*D\rbrack$, with
 $V^*_i$ replaced by $V^*_{d-i}$
 for
$0 \leq i \leq d$.
\hfill $\Box $ \\

\section{The linear transformations $B, B^*, K, K^*$}

\noindent In the previous two sections we discussed
general tridiagonal pairs. For the rest of this
paper we restrict our attention to the special case
mentioned in Theorem
\ref{thm:main}. We will refer to the following setup.

\begin{definition}
\label{def:setup2}
\rm
Let $V$ denote a vector space over $\K$ with finite positive
dimension  and let $A,A^*$ denote a tridiagonal pair on $V$.
Let $V_0, V_1, \ldots, V_d$ (resp.  
$V^*_0, V^*_1, \ldots, V^*_d$)
denote a standard ordering of the eigenspaces of $A$
(resp. $A^*$).
For $0 \leq i\leq d$ let $\theta_i$ (resp. $\theta^*_i$)
denote the eigenvalue of $A$ (resp. $A^*$) associated
with $V_i$ (resp. $V^*_i$).
We assume
there exist nonzero scalars $a, a^*$ in $\K$ such that
\begin{eqnarray}
\theta_i = aq^{2i-d}, \qquad \qquad 
\theta^*_i = a^*q^{d-2i} \qquad \qquad 
(0\leq i\leq d).
\label{eq:eig}
\end{eqnarray}
Let $b$ and $b^*$ denote nonzero scalars in $\K$.
\end{definition}

\begin{definition}
\label{def:b}
\label{def:k}
\rm
Adopt the assumptions of Definition \ref{def:setup2}.
\begin{enumerate}
\item 
We let $B:V\rightarrow V$ denote the unique
linear transformation
such that
for $0 \leq i \leq d$, 
\begin{eqnarray}
(V^*_0+\cdots+ V^*_i)\cap (V_0+\cdots + V_{d-i})
\label{eq:oos}
\end{eqnarray}
is an eigenspace of $B$ with eigenvalue 
$bq^{2i-d}$. We remark
(\ref{eq:oos}) is the $i$th subspace of the decomposition
$\lbrack 0^* 0\rbrack $ from Lemma 
\ref{thm:sixdecp}.
\item
We let $B^*:V\rightarrow V$ denote the unique
linear transformation
such that
for $0 \leq i \leq d$, 
\begin{eqnarray}
(V^*_{d-i}+\cdots +V^*_d)\cap (V_i+\cdots + V_d)
\label{eq:oos2}
\end{eqnarray}
is an eigenspace of $B^*$ with eigenvalue 
$b^*q^{d-2i}$. We remark
(\ref{eq:oos2}) is the $i$th subspace of the decomposition
$\lbrack D^* D\rbrack $ from Lemma 
\ref{thm:sixdecp}.
\item
We let $K:V\rightarrow V$ denote the unique
linear transformation
such that
for $0 \leq i \leq d$, 
\begin{eqnarray}
(V^*_0+\cdots +V^*_i)\cap (V_i+\cdots + V_d)
\label{eq:ksp}
\end{eqnarray}
is an eigenspace of $K$ with eigenvalue 
$q^{2i-d}$. We remark
(\ref{eq:ksp}) is the $i$th subspace of the decomposition
$\lbrack 0^* D\rbrack $ from Lemma 
\ref{thm:sixdecp}.
\item
 We let $K^*:V\rightarrow V$ denote the unique
linear transformation
such that
for $0 \leq i \leq d$, 
\begin{eqnarray}
(V^*_{d-i}+\cdots +V^*_d)\cap (V_0+\cdots + V_{d-i})
\label{eq:kspd}
\end{eqnarray}
is an eigenspace of $K^*$ with eigenvalue 
$q^{2i-d}$. We remark
(\ref{eq:kspd}) is the $i$th subspace of the decomposition
$\lbrack D^* 0\rbrack $ from Lemma 
\ref{thm:sixdecp}.
\end{enumerate}
\end{definition}

\begin{remark}
\rm
\label{rem:inv}
With reference to
Definition
 \ref{def:setup2} and Definition
\ref{def:b}, the following (i), (ii) hold.
\begin{enumerate}
\item
If we replace 
$(A,A^*,V_i, V^*_i, $ $a,a^*,B,B^*, b,b^*,K,K^*,q)$
by 
$(A^*,A,V^*_{d-i}, V_{d-i}, a^*,a,
B^*,B,$ $ b^*,b,K^{-1}, K^{*-1},q)$ 
then the requirements of
 Definition
 \ref{def:setup2} and Definition
\ref{def:b}
 are still satisfied.
\item
If we replace 
$(A,A^*,V_i, V^*_i, $ $a,a^*,B,B^*, b,b^*,K,K^*,q)$
by 
$(A,A^*,V_{d-i}, V^*_{d-i}, a,a^*,
B^*,B, $ $ b^*,b,K^{*-1}, K^{-1}, q^{-1})$
then the requirements of
 Definition
 \ref{def:setup2} and Definition
\ref{def:b}
 are still satisfied.
\end{enumerate}
\end{remark}
We will use
Remark
\ref{rem:inv}
to streamline a few proofs later in the paper.

\section{Some relations involving $A, A^*,B,B^*$}

\noindent In this section we give four relations
involving the tridiagonal pair $A,A^*$ 
from Definition
\ref{def:setup2} and the elements $B,B^*$ from
Definition
\ref{def:b}.

\begin{theorem}
\label{thm:key}
With reference to Definition \ref{def:setup2}
and Definition
\ref{def:b},
\begin{eqnarray}
\frac{qAB-q^{-1}BA}{q-q^{-1}}&=&abI,
\label{eq:eq1}
\\
\frac{qBA^*-q^{-1}A^*B}{q-q^{-1}}&=&a^*bI,
\label{eq:eq2}
\\
\frac{qA^*B^*-q^{-1}B^*A^*}{q-q^{-1}}&=&a^*b^*I,
\label{eq:eq4}
\\
\frac{qB^*A-q^{-1}AB^*}{q-q^{-1}}&=&ab^*I.
\label{eq:eq3}
\end{eqnarray}
\end{theorem}
\noindent {\it Proof:}
We first show 
(\ref{eq:eq1}).
Let $U_0, U_1, \ldots, U_d$ denote the decomposition
$\lbrack 0^*0\rbrack $ from
Lemma
\ref{thm:sixdecp}.
We show
$qAB-q^{-1}BA-ab(q-q^{-1})I$ 
vanishes on $U_i$ for
$0 \leq i \leq d$.
Let $i$ be given.
By Definition
\ref{def:b}(i) we find $B-bq^{2i-d}I$ vanishes on
 $U_i$ so
\begin{eqnarray}
\label{eq:p1}
(A-aq^{d-2i-2}I)(B-bq^{2i-d}I)
\end{eqnarray}
vanishes on
 $U_i$.
From the table of Lemma
\ref{thm:aaction}, row $\lbrack 0^*0\rbrack $,
and using 
(\ref{eq:eig}),
we find
$(A-aq^{d-2i}I)U_i\subseteq U_{i+1}$. Therefore
\begin{eqnarray}
\label{eq:p2}
(B-bq^{2i+2-d}I)(A-aq^{d-2i}I)
\end{eqnarray}
vanishes on $U_i$.
Subtracting $q^{-1}$times 
(\ref{eq:p2}) from
$q$ times (\ref{eq:p1}) 
we find
$qAB-q^{-1}BA-ab(q-q^{-1})I$ vanishes on
 $U_i$.
 Line
(\ref{eq:eq1}) follows.
To get 
(\ref{eq:eq4}) use
(\ref{eq:eq1}) 
and the involution given in 
Remark \ref{rem:inv}(i).
To get 
(\ref{eq:eq3}) use
(\ref{eq:eq1}) 
and the involution
given in
Remark \ref{rem:inv}(ii).
To get 
(\ref{eq:eq2}) use
(\ref{eq:eq4}) 
and the involution given in 
Remark \ref{rem:inv}(ii).
\hfill $\Box $ \\

\section{The action of $B$ and $B^*$ on the
six decompositions}

\noindent In this section we describe how the elements
$B,B^*$ from Definition \ref{def:b} act on 
the six decompositions given in Lemma
\ref{thm:sixdecp}.

\begin{theorem}
\label{thm:baction}
Adopt the assumptions of Definition \ref{def:setup2}
and let $U_0, U_1, \ldots, U_d$ 
denote any one of the six  decompositions
of $V$ given in 
Lemma \ref{thm:sixdecp}.
Let the maps $B, B^*$ be as in  Definition \ref{def:b}.
Then for $0 \leq i \leq d$
the action of $B$ and $B^*$ on $U_i$ is described as follows.

\medskip

\centerline{
\begin{tabular}[t]{c|c|c}
      {\rm name} &{\rm action of $B$ on $U_i$} & {\rm action of $B^*$ on $U_i$}
      \\ \hline  \hline
	$\lbrack 0D\rbrack $ &
	$ (B-bq^{d-2i}I)U_i\subseteq U_{i-1}$  & 
	$(B^*-b^*q^{d-2i}I)U_i \subseteq U_{i+1}$ 
\\	
	$\lbrack 0^*D^*\rbrack $ &
	$(B-bq^{2i-d}I)U_i \subseteq U_{i-1}$  &
	$ (B^*-b^*q^{2i-d}I)U_i\subseteq U_{i+1}$   \\
       $\lbrack 0^*D\rbrack $ &
        $(B-bq^{2i-d}I)U_i\subseteq U_{i-1}$ &
	$(B^*-b^*q^{d-2i}I)U_i \subseteq U_{i+1}$   \\ 
	$\lbrack 0^*0\rbrack $ & 
	$(B-bq^{2i-d}I)U_i=0$ &
	$B^*U_i \subseteq U_{i-1}+U_i + U_{i+1}$   \\ 
       $\lbrack D^*0\rbrack $ & 
	$(B-bq^{2i-d}I)U_i\subseteq U_{i+1}$ &
	$(B^*-b^*q^{d-2i}I)U_i \subseteq U_{i-1}$   \\ 
	$\lbrack D^*D\rbrack $ &
	$BU_i\subseteq U_{i-1}+U_i + U_{i+1}$ &
	$(B^*-b^*q^{d-2i}I)U_i=0$  
	\end{tabular}}
\medskip

\end{theorem}
\noindent {\it Proof:}
We first give the action of $B$ for each of the six
rows in the table.
\\
\noindent 
$\lbrack 0D\rbrack $: 
Let $U_0, U_1, \ldots, U_d$ denote the decomposition
$\lbrack 0D\rbrack $. 
From Lemma
\ref{thm:aaction},
row $\lbrack 0D\rbrack $, and using
(\ref{eq:eig}),
we find that for $0 \leq i \leq d$,  
$U_i$ is an eigenspace for $A$ with
eigenvalue $aq^{2i-d}$.
We show
$(B-bq^{d-2i}I)U_i \subseteq U_{i-1}$ for  $0 \leq i \leq d$.
To do this, it suffices to show
\begin{eqnarray}
(A-aq^{2i-2-d}I)(B-bq^{d-2i}I)
\label{eq:want}
\end{eqnarray}
vanishes on $U_i$ for $0 \leq i \leq d$.
Let $i$ be given.
Observe
 $A-aq^{2i-d}I$ vanishes on
$U_i$ so
\begin{eqnarray}
(B-bq^{d-2i+2}I)(A-aq^{2i-d}I)
\label{eq:ba}
\end{eqnarray}
vanishes on $U_i$. Using 
(\ref{eq:eq1}) we find
\begin{eqnarray}
qAB-q^{-1}BA-ab(q-q^{-1})I
\label{eq:qab}
\end{eqnarray}
vanishes on $U_i$.
Adding
(\ref{eq:ba})
to 
$q$ times 
(\ref{eq:qab}) 
we find
(\ref{eq:want}) vanishes on $U_i$.
We conclude $(B-bq^{d-2i}I)U_i \subseteq U_{i-1}$ for  $0 \leq i \leq d$.
\\
$\lbrack 0^*D^*\rbrack $:
Let $U_0, U_1, \ldots, U_d$ denote the decomposition
$\lbrack 0^*D^*\rbrack $.
From Lemma
\ref{thm:aaction},
row  
$\lbrack 0^*D^*\rbrack $, and using
(\ref{eq:eig}), 
we find that for $0 \leq i \leq d$,  
$U_i$ is an eigenspace for $A^*$ with
eigenvalue $a^*q^{d-2i}$.
We show
$(B-bq^{2i-d}I)U_i \subseteq U_{i-1}$ for $0 \leq i \leq d$.
To do this,  it suffices to show
\begin{eqnarray}
(A^*-a^*q^{d-2i+2}I)(B-bq^{2i-d}I)
\label{eq:wantdual}
\end{eqnarray}
vanishes on $U_i$ for $0 \leq i\leq d$. 
Let $i$ be given.
Observe
$A^*-a^*q^{d-2i}I$ vanishes on
$U_i$ so
\begin{eqnarray}
(B-bq^{2i-d-2}I)(A^*-a^*q^{d-2i}I)
\label{eq:badual}
\end{eqnarray}
vanishes on $U_i$. Using 
(\ref{eq:eq2}) we find
\begin{eqnarray}
qBA^*-q^{-1}A^*B-a^*b(q-q^{-1})I
\label{eq:qabdual}
\end{eqnarray}
vanishes on $U_i$.
Subtracting (\ref{eq:badual})
from
$q^{-1}$ times 
(\ref{eq:qabdual})
we find
(\ref{eq:wantdual}) vanishes on $U_i$.
We conclude 
$(B-bq^{2i-d}I)U_i \subseteq U_{i-1}$ for $0 \leq i \leq d$.
\\
\noindent $\lbrack 0^*D\rbrack $:
Let $U_0, U_1, \ldots, U_d$ denote the decomposition
 $\lbrack 0^*D\rbrack $.
We show
$(B-bq^{2i-d}I)U_i \subseteq U_{i-1}$ for $0 \leq i \leq d$.
Let $i$ be given. We have
\begin{eqnarray*}
(B-bq^{2i-d}I)U_i  &\subseteq &
(B-bq^{2i-d}I)(U_0+\cdots +U_{i})
\\
&=& (B-bq^{2i-d}I)(V^*_0+\cdots + V^*_i)
\qquad  (\mbox{by
Lemma 
\ref{thm:decsum}, row $\lbrack 0^*D\rbrack$})
\\
&\subseteq& V^*_0+\cdots +V^*_{i-1}
\qquad \qquad \qquad (\mbox{by
present Theorem, row
 $\lbrack 0^*D^*\rbrack $})
\\
&=& U_0+\cdots +U_{i-1}
\qquad \qquad
\qquad \qquad (\mbox{by
Lemma 
\ref{thm:decsum}, row $\lbrack 0^*D\rbrack $})
\end{eqnarray*}
and also
\begin{eqnarray*}
(B-bq^{2i-d}I)U_i  &\subseteq &
(B-bq^{2i-d}I)(U_i+\cdots + U_{d})
\\
&=& (B-bq^{2i-d}I)(V_i+\cdots +V_d)
\qquad  (\mbox{by
Lemma 
\ref{thm:decsum}, row $\lbrack 0^*D\rbrack $})
\\
&\subseteq & V_{i-1}+\cdots +V_d
\qquad \qquad \qquad (\mbox{by
present Theorem, row
 $\lbrack 0D\rbrack $})
\\
&=& U_{i-1}+\cdots +U_d
\qquad \qquad 
\qquad \qquad (\mbox{by
Lemma 
\ref{thm:decsum}, row $\lbrack 0^*D\rbrack $}).
\end{eqnarray*}
Combining these observations 
we obtain
$(B-bq^{2i-d}I)U_i \subseteq U_{i-1}$ for $0 \leq i \leq d$.
\\
\noindent 
 $\lbrack 0^*0\rbrack $:
Let 
$U_0, U_1, \ldots, U_d$ denote the decomposition
 $\lbrack 0^*0\rbrack $.
Then $(B-bq^{2i-d}I)U_i=0$ for $0 \leq i \leq d$
by Definition
\ref{def:b}(i).
\\
\noindent
 $\lbrack D^*0\rbrack $:
Let $U_0, U_1, \ldots, U_d$ denote the decomposition
 $\lbrack D^*0\rbrack $.
We show
$(B-bq^{2i-d}I)U_i \subseteq U_{i+1}$ for $0 \leq i \leq d$.
Let $i$ be given. We have
\begin{eqnarray*}
(B-bq^{2i-d}I)U_i  &\subseteq &
(B-bq^{2i-d}I)(U_0+\cdots +U_{i})
\\
&=& (B-bq^{2i-d}I)(V^*_{d-i}+\cdots + V^*_d)
\qquad (\mbox{by
Lemma 
\ref{thm:decsum}, row $\lbrack D^*0\rbrack $})
\\
&\subseteq& V^*_{d-i-1}+\cdots +V^*_d
\qquad \qquad \qquad (\mbox{by
present Theorem, row
 $\lbrack 0^*D^*\rbrack $})
\\
&=& U_0+\cdots +U_{i+1}
\qquad \qquad  \quad 
\qquad \qquad (\mbox{by
Lemma 
\ref{thm:decsum}, row $\lbrack D^*0\rbrack $})
\end{eqnarray*}
and also
\begin{eqnarray*}
(B-bq^{2i-d}I)U_i  &\subseteq &
(B-bq^{2i-d}I)(U_i+\cdots + U_{d})
\\
&=& (B-bq^{2i-d}I)(V_0+\cdots +V_{d-i})
\qquad (\mbox{by
Lemma 
\ref{thm:decsum}, row $\lbrack D^*0\rbrack $})
\\
&\subseteq & V_{0}+\cdots +V_{d-i-1}
\qquad \qquad \qquad (\mbox{by
present Theorem, row
 $\lbrack 0D\rbrack $})
\\
&=& U_{i+1}+\cdots +U_{d}
\qquad \qquad  \quad 
\qquad \qquad (\mbox{by
Lemma 
\ref{thm:decsum}, row $\lbrack D^*0\rbrack $}).
\end{eqnarray*}
Combining these observations 
we obtain
$(B-bq^{2i-d}I)U_i \subseteq U_{i+1}$ for $0 \leq i \leq d$.
\\
\noindent 
 $\lbrack D^*D\rbrack $:
Let 
$U_0, U_1, \ldots, U_d$ denote the decomposition
 $\lbrack D^*D\rbrack $.
We show
$B U_i \subseteq 
 U_{i-1} 
+
 U_i 
+
 U_{i+1}$
 for $0 \leq i \leq d$.
Let $i$ be given. We have
\begin{eqnarray*}
B U_i &\subseteq&
B (U_0+ \cdots + U_i)
\\
&=& B (V^*_{d-i}+ \cdots + V^*_d)
\qquad \qquad (\mbox{by 
Lemma
\ref{thm:decsum}, row $\lbrack D^*D\rbrack $})
\\
&\subseteq& V^*_{d-i-1}+ \cdots + V^*_d
\qquad \qquad  (\mbox{by 
present Theorem,
row $\lbrack 0^*D^*\rbrack $})
\\
&=& 
U_0+ \cdots + U_{i+1}
\qquad \qquad 
\qquad \qquad (\mbox{by 
Lemma
\ref{thm:decsum}, row $\lbrack D^*D\rbrack $})
\end{eqnarray*}
and also
\begin{eqnarray*}
B U_i &\subseteq&
B (U_i+ \cdots + U_d)
\\
&=& B (V_i+ \cdots + V_d)
\qquad \qquad (\mbox{by 
Lemma
\ref{thm:decsum}, row $\lbrack D^*D\rbrack $})
\\
&\subseteq& V_{i-1}+ \cdots + V_d
\qquad \qquad  (\mbox{by 
present Theorem, 
row $\lbrack 0D\rbrack $})
\\
&=& 
U_{i-1}+ \cdots + U_d
\qquad 
\qquad \qquad (\mbox{by 
Lemma
\ref{thm:decsum}, row $\lbrack D^*D\rbrack $}).
\end{eqnarray*}
Combining these observations 
we 
find
$B U_i \subseteq 
 U_{i-1} 
+
 U_i 
+
 U_{i+1}$ for $0 \leq i \leq d$.
\\
\noindent We have now given the action of $B$ on each
of the six decompositions.
Using this and the involution 
from Remark
\ref{rem:inv}(i),
we find $B^*$ acts on
the six  decompositions as claimed.
\hfill $\Box $ \\

\section{The pair $B,B^*$ is a tridiagonal pair}

\noindent In this section we show 
that the linear transformations $B,B^*$ from Definition
\ref{def:b} form a tridiagonal pair.

\begin{theorem}
Adopt the assumptions of Definition \ref{def:setup2}
and let the maps $B,B^*$ be as in Definition \ref{def:b}.
Then the pair $B,B^*$ is a tridiagonal pair on $V$.
The sequence $bq^{2i-d}$ $(0 \leq i \leq d)$
is a standard ordering of the eigenvalues of $B$
and 
the sequence $b^*q^{d-2i}$ $(0 \leq i \leq d)$
is a standard ordering of the eigenvalues of $B^*$.
\end{theorem}
\noindent {\it Proof:} For the duration of this proof
 let $U_0,  \ldots, U_d$ 
(resp. $U^*_0, \ldots, U^*_d$)
denote the decomposition $\lbrack 0^*0\rbrack$ 
(resp. $\lbrack D^*D\rbrack$) from 
Lemma 
\ref{thm:sixdecp}.
We show the pair $B,B^*$ is a tridiagonal pair on $V$.
To do this we show $B,B^*$ satisfies conditions (i)--(iv)
in Definition 
\ref{def:tdp}.
\\
\noindent {\it Proof that $B,B^*$ satisfies 
Definition 
\ref{def:tdp}(i)}:
Each of $U_0, \ldots,U_d$ is an eigenspace of
$B$ by 
 Definition \ref{def:b}(i) and these eigenspaces span $V$  
so $B$ is diagonalizable.
Each of $U^*_0, \ldots,U^*_d$ is an eigenspace of
$B^*$ by 
 Definition \ref{def:b}(ii) and these eigenspaces span $V$ 
so $B^*$ is diagonalizable.
\\
\noindent {\it Proof that $B,B^*$ satisfies 
Definition 
\ref{def:tdp}(ii)}:
From the construction
$U_0, \ldots, U_d$ is an ordering of the eigenspaces
of $B$. By
Theorem
\ref{thm:baction}, row $\lbrack 0^*0\rbrack $ we find
$B^*U_i \subseteq U_{i-1}+U_i + U_{i+1}$
for $0 \leq i \leq d$.
\\
\noindent {\it Proof that $B,B^*$ satisfies 
Definition 
\ref{def:tdp}(iii)}:
From the construction
$U^*_0, \ldots, U^*_d$ is an ordering of the eigenspaces
of $B^*$. By
Theorem
\ref{thm:baction}, row $\lbrack D^*D\rbrack $ we find
$BU^*_i \subseteq U^*_{i-1}+U^*_i + U^*_{i+1}$
for $0 \leq i \leq d$.
\\
\noindent {\it Proof that $B,B^*$ satisfies 
Definition 
\ref{def:tdp}(iv)}:
We let $W$ denote an irreducible
$(B,B^*)$-submodule of $V$ and show
$W=V$. To obtain $W=V$ we will show $A W\subseteq W$ and
$A^*W\subseteq W$.
We first show $AW\subseteq W$. 
We define ${\tilde W}:=\lbrace w \in W | Aw \in W\rbrace $
and show ${\tilde W}=W$.
Using
(\ref{eq:eq1}) we routinely find $B{\tilde W}\subseteq {\tilde W}$.
Using
(\ref{eq:eq3}) we routinely find $B^*{\tilde W}\subseteq {\tilde W}$.
We claim  
$\tilde W\not=0$.
To prove the claim,
We define
$W_i = W \cap U_i$
for $0 \leq i \leq d$.
From the table of Lemma
\ref{thm:decsum}, row $\lbrack 0^*0\rbrack $
we find both
\begin{eqnarray}
\label{eq:s1}
W_0+ \cdots + W_i \subseteq V^*_0 + \cdots + V^*_i
&\qquad & (0 \leq i \leq d),
\\
\label{eq:s2}
W_i+ \cdots + W_d \subseteq V_0 + \cdots + V_{d-i}
&\qquad & (0 \leq i \leq d).
\end{eqnarray}
The nonzero spaces among $W_0, \ldots, W_d$ are the eigenspaces 
of $B$ on $W$ so
$W=\sum_{i=0}^d W_i$.
By this and since $W\not=0$ we find
$W_0, \ldots, W_d$ are not all 0.
Define $r=\mbox{max}\lbrace i |0 \leq i \leq d, \;W_i\not=0\rbrace$.
We define 
$W^*_i = W\cap U^*_i$
for $0 \leq i \leq d$.
From the table of
Lemma \ref{thm:decsum}, row $\lbrack D^*D\rbrack $ we find
\begin{eqnarray}
\label{eq:s3}
W^*_0+ \cdots + W^*_i \subseteq V^*_{d-i}+ \cdots + V^*_d
&\qquad & (0 \leq i \leq d),
\\
\label{eq:s4}
W^*_i+ \cdots + W^*_d \subseteq V_i + \cdots + V_d
&\qquad & (0 \leq i \leq d).
\end{eqnarray}
The nonzero spaces among $W^*_0, \ldots, W^*_d$ are the eigenspaces 
of $B^*$ on $W$ so
$W=\sum_{i=0}^d W^*_i$.
By this and since $W\not=0$
we find
$W^*_0, \ldots, W^*_d$ are not all 0.
Define $t=\mbox{min}\lbrace i |0 \leq i \leq d, \;W^*_i\not=0\rbrace$.
Suppose for the moment that $r+t<d$.
Setting $i=r$ in
(\ref{eq:s1}) 
and using $W_0+\cdots + W_r=W$ 
we find
$W\subseteq V^*_0+\cdots +V^*_r$.
Setting $i=t$ in 
(\ref{eq:s3}) we find $W^*_t\subseteq V^*_{d-t}+\cdots +V^*_d$.
Of course $W^*_t \subseteq W$ so 
\begin{eqnarray*}
W^*_t &=& W \cap W^*_t
\\
&\subseteq&
(V^*_0+\cdots +V^*_r)\cap 
(V^*_{d-t}+\cdots +V^*_d)
\\
&=& 0
\end{eqnarray*}
for a contradiction. Therefore $r+t\geq d$.
Setting $i=r$ in
(\ref{eq:s2}) we find $W_r\subseteq V_0+\cdots +V_{d-r}$.
Setting $i=t$ in
(\ref{eq:s4}) and using 
$W^*_t + \cdots + W^*_d=W$ 
we find
$W \subseteq V_t + \cdots + V_d$.
Of course $W_r \subseteq W$ so
\begin{eqnarray*}
W_r &=& W_r\cap W
\\
&\subseteq &
(V_0+\cdots +V_{d-r})\cap 
(V_t+\cdots +V_d).
\end{eqnarray*}
By this and since $r+t\geq d$ we 
find $r+t=d$ and then $W_r\subseteq V_{d-r}$.
Recall $V_{d-r}$ is an eigenspace for
$A$ so 
$AW_r\subseteq W_r$.
Therefore
$AW_r\subseteq W$ so
$W_r \subseteq {\tilde W}$.
Consequently ${\tilde W}\not=0$ as desired.
We have shown ${\tilde W}$ is nonzero and invariant under each of
$B,B^*$.
Therefore ${\tilde W}=W$
since  $W$ is irreducible  as a $(B,B^*)$-module.
We have now shown $AW\subseteq W$.
Using this and the involution in Remark
\ref{rem:inv}(i) 
we find $A^*W\subseteq W$.
Applying
Definition 
\ref{def:tdp}(iv) to $A,A^*$ we find $W=V$.
\\
\noindent We have now shown
the pair $B,B^*$ satisfies 
conditions (i)--(iv) of 
Definition \ref{def:tdp}.
Therefore $B,B^*$ is a tridiagonal pair on $V$.
From the construction 
$U_0,\ldots, U_d$ is a standard ordering
of the eigenspaces of $B$.
For $0 \leq i \leq d$ the scalar
 $bq^{2i-d}$ is the eigenvalue of $B$ associated with $U_i$.
Therefore the sequence $bq^{2i-d}$ $(0\leq i \leq d)$
is a standard ordering of the eigenvalues of $B$.
From the construction 
$U^*_0,\ldots, U^*_d$ is a standard ordering
of the eigenspaces of $B^*$.
For $0 \leq i \leq d$ the scalar
 $b^*q^{d-2i}$ is the eigenvalue of $B^*$ associated with $U^*_i$.
Therefore the sequence $b^*q^{d-2i}$ $(0 \leq i\leq d)$
is a standard ordering of the eigenvalues of $B^*$.
\hfill $\Box $ \\

\section{Some relations involving $A, A^*,B,B^*, K,K^*$}

\noindent In this section we give some relations involving
the tridiagonal pair $A,A^*$ from
Definition
\ref{def:setup2}, the tridiagonal pair $B,B^*$ from
Definition \ref{def:b}, and the elements
$K, K^*$ from Definition
\ref{def:k}.

\begin{theorem}
\label{thm:kkey}
With reference to Definition \ref{def:setup2} and
Definition \ref{def:b}, 
\begin{eqnarray}
\frac{qK^{-1}A-q^{-1}AK^{-1}}{q-q^{-1}}&=&aI,
\label{eq:keq1}
\\
\frac{qBK^{-1}-q^{-1}K^{-1}B}{q-q^{-1}}&=&bI,
\label{eq:keq2}
\\
\frac{qKA^*-q^{-1}A^*K}{q-q^{-1}}&=&a^*I,
\label{eq:keq4}
\\
\frac{qB^*K-q^{-1}KB^*}{q-q^{-1}}&=&b^*I.
\label{eq:keq3}
\end{eqnarray}
\end{theorem}
\noindent {\it Proof:}
We first show 
(\ref{eq:keq1}),
(\ref{eq:keq2}).
Let $U_0, U_1, \ldots, U_d$ denote the decomposition
$\lbrack 0^*D\rbrack $ from
Lemma
\ref{thm:sixdecp}.
Concerning
(\ref{eq:keq1}), 
we show
$qK^{-1}A-q^{-1}AK^{-1}-a(q-q^{-1})I$ 
vanishes on $U_i$ for
$0 \leq i \leq d$.
Let $i$ be given.
Observe $K-q^{2i-d}I$ vanishes on
 $U_i$ by 
 Definition
\ref{def:k} so $K^{-1}-q^{d-2i}I$ vanishes on
 $U_i$; from this we find
\begin{eqnarray}
\label{eq:kp1}
(A-aq^{2i-d+2}I)(K^{-1}-q^{d-2i}I)
\end{eqnarray}
vanishes on
 $U_i$.
From the table of Lemma
\ref{thm:aaction}, row $\lbrack 0^*D\rbrack $,
and using (\ref{eq:eig}),
we find
$(A-aq^{2i-d}I)U_i\subseteq U_{i+1}$. Therefore
\begin{eqnarray}
\label{eq:kp2}
(K^{-1}-q^{d-2i-2}I)(A-aq^{2i-d}I)
\end{eqnarray}
vanishes on $U_i$.
Subtracting $q^{-1}$times 
(\ref{eq:kp1}) from
$q$ times (\ref{eq:kp2}) 
we find
$qK^{-1}A-q^{-1}AK^{-1}-a(q-q^{-1})I$ 
 vanishes on $U_i$.
 Line
(\ref{eq:keq1}) follows.
Concerning
(\ref{eq:keq2}), 
we show
$qBK^{-1}-q^{-1}K^{-1}B-b(q-q^{-1})I$ vanishes on $U_i$ for
$0 \leq i \leq d$.
Let $i$ be given.
We mentioned earlier that $K^{-1}-q^{d-2i}I$ vanishes on
 $U_i$
so
\begin{eqnarray}
\label{eq:kp1n}
(B-bq^{2i-d-2}I)(K^{-1}-q^{d-2i}I)
\end{eqnarray}
vanishes on
 $U_i$.
From the table of Lemma
\ref{thm:baction}, row $\lbrack 0^*D\rbrack $,
we find
$(B-bq^{2i-d}I)U_i \subseteq U_{i-1}$. Therefore 
\begin{eqnarray}
\label{eq:kp2n}
(K^{-1}-q^{d-2i+2}I)(B-bq^{2i-d}I)
\end{eqnarray}
vanishes on $U_i$.
Subtracting $q^{-1}$times 
(\ref{eq:kp2n}) from
$q$ times (\ref{eq:kp1n}) 
we find
$qBK^{-1}-q^{-1}K^{-1}B-b(q-q^{-1})I$ vanishes on $U_i$.
Line
(\ref{eq:keq2}) follows.
To obtain  
(\ref{eq:keq4}), 
(\ref{eq:keq3})
apply
(\ref{eq:keq1}), (\ref{eq:keq2})  and the involution
given in Remark \ref{rem:inv}(i). 
\hfill $\Box $ \\

\begin{theorem}
\label{thm:kkeyd}
With reference to Definition \ref{def:setup2} and
Definition \ref{def:b}, 
\begin{eqnarray}
\frac{qAK^*-q^{-1}K^*A}{q-q^{-1}}&=&aI,
\label{eq:keq3d}
\\
\frac{qK^{*-1}B-q^{-1}BK^{*-1}}{q-q^{-1}}&=&bI,
\label{eq:keq1d}
\\
\frac{qA^*K^{*-1}-q^{-1}K^{*-1}A^*}{q-q^{-1}}&=&a^*I,
\label{eq:keq2d}
\\
\frac{qK^{*}B^*-q^{-1}B^*K^*}{q-q^{-1}}&=&b^*I.
\label{eq:keq4d}
\end{eqnarray}
\end{theorem}
\noindent {\it Proof:}
Use Theorem
\ref{thm:kkey} and the involution
given in
Remark
\ref{rem:inv}(ii).
\hfill $\Box $ \\

\section{ The actions of $K,K^{-1}, K^*, K^{*-1}$ on
the six decompositions}

\noindent In this section  we describe how 
the elements $K,K^{-1}, K^*, K^{*-1}$ from
Definition
\ref{def:k}
act on the six decompositions
from
Lemma
\ref{thm:sixdecp}.
We begin with $K$ and $K^{-1}$.

\begin{theorem}
\label{thm:kkdaction}
Adopt the assumptions of Definition \ref{def:setup2} and
let $U_0, U_1, \ldots, U_d$ 
denote any one of the six decompositions
of $V$ given in 
Lemma \ref{thm:sixdecp}. Let the map $K$
be as in 
Definition \ref{def:b}.
Then for $0 \leq i \leq d$
the action of $K$ and $K^{-1}$ on $U_i$ is described as follows.

\medskip
\centerline{
\begin{tabular}[t]{c|c|c}
      {\rm name} &{\rm action of $K$ on $U_i$} & {\rm action of $K^{-1}$ on $U_i$}
      \\ \hline  \hline
	$\lbrack 0D\rbrack $ &
	$ (K-q^{2i-d}I)U_i\subseteq U_{i+1}+\cdots + U_d$  & 
	$(K^{-1}-q^{d-2i}I) U_i \subseteq U_{i+1}$ 
\\	
	$\lbrack 0^*D^*\rbrack $ &
	$(K- q^{2i-d}I)U_i \subseteq U_{i-1}$  &
	$ (K^{-1}-q^{d-2i}I)U_i\subseteq U_0+\cdots + U_{i-1}$   \\
       $\lbrack 0^*D\rbrack $ &
        $(K-q^{2i-d}I)U_i=0$ &
	$(K^{-1}-q^{d-2i}I)U_i=0$   \\ 
	$\lbrack 0^*0\rbrack $ & 
	$(K-q^{2i-d}I)U_i\subseteq U_0+\cdots +U_{i-1}$ &
	$(K^{-1}-q^{d-2i}I)U_i \subseteq U_{i-1}$   \\ 
       $\lbrack D^*0\rbrack $ & 
	$KU_i\subseteq U_0+\cdots+ U_{i+1}$ &
	$K^{-1}U_i \subseteq U_{i-1}+\cdots + U_d$   \\ 
	$\lbrack D^*D\rbrack $ &
	$(K-q^{2i-d}I)U_i\subseteq U_{i+1}$ &
	$(K^{-1}-q^{d-2i}I)U_i \subseteq U_{i+1}+\cdots+U_d$  
	\end{tabular}}
\medskip
\end{theorem}
\noindent {\it Proof:}
We consider each of the six
rows of the table.
\\
\noindent 
$\lbrack 0D\rbrack $: 
Let $U_0, U_1, \ldots, U_d$ denote the decomposition
$\lbrack 0D\rbrack $. 
From Lemma
\ref{thm:aaction},
row $\lbrack 0D\rbrack $, and using
(\ref{eq:eig}),
we find that for $0 \leq i \leq d$,  
$U_i$ is an eigenspace for $A$ with
eigenvalue $aq^{2i-d}$.
We show
$(K^{-1}-q^{d-2i}I)U_i \subseteq U_{i+1}$ for  $0 \leq i \leq d$.
To do this, it suffices to show
\begin{eqnarray}
(A-aq^{2i+2-d}I)(K^{-1}-q^{d-2i}I)
\label{eq:wantk}
\end{eqnarray}
vanishes on $U_i$ for $0 \leq i \leq d$.
Let $i$ be given.
Observe
 $A-aq^{2i-d}I$ vanishes on
$U_i$ so
\begin{eqnarray}
(K^{-1}-q^{d-2i-2}I)(A-aq^{2i-d}I)
\label{eq:bak}
\end{eqnarray}
vanishes on $U_i$. Using 
(\ref{eq:keq1}) we find
\begin{eqnarray}
qK^{-1}A-q^{-1}AK^{-1}-a(q-q^{-1})I
\label{eq:qabk}
\end{eqnarray}
vanishes on $U_i$.
Subtracting
(\ref{eq:bak})
from
$q^{-1}$ times 
(\ref{eq:qabk}) 
we find
(\ref{eq:wantk}) vanishes on $U_i$.
We conclude 
$(K^{-1}-q^{d-2i}I)U_i \subseteq U_{i+1}$ for  $0 \leq i \leq d$.
From this
we find
$(K-q^{2i-d}I)U_i \subseteq U_{i+1}+\cdots + U_d$ for  $0 \leq i \leq d$.
\\
$\lbrack 0^*D^*\rbrack $:
Use the present Theorem, row 
$\lbrack 0D\rbrack $ and the involution
given in Remark \ref{rem:inv}(i).
\\
\noindent $\lbrack 0^*D\rbrack $:
Let $U_0, U_1, \ldots, U_d$ denote the decomposition
 $\lbrack 0^*D\rbrack $.
From Definition \ref{def:k} we find
$(K-q^{2i-d}I)U_i =0$ for $0 \leq i \leq d$.
It follows
$(K^{-1}-q^{d-2i}I)U_i =0$ for $0 \leq i \leq d$.
\\
\noindent 
 $\lbrack 0^*0\rbrack $:
Let 
$U_0, U_1, \ldots, U_d$ denote the decomposition
 $\lbrack 0^*0\rbrack $.
From 
Definition \ref{def:b}
we find that for $0 \leq i \leq d$,  
$U_i$ is an eigenspace for $B$ with
eigenvalue $bq^{2i-d}$.
We show
$(K^{-1}-q^{d-2i}I)U_i \subseteq U_{i-1}$ for  $0 \leq i \leq d$.
To do this, it suffices to show
\begin{eqnarray}
(B-bq^{2i-2-d}I)(K^{-1}-q^{d-2i}I)
\label{eq:wantkb}
\end{eqnarray}
vanishes on $U_i$ for $0 \leq i \leq d$.
Let $i$ be given.
Observe
 $B-bq^{2i-d}I$ vanishes on
$U_i$ so
\begin{eqnarray}
(K^{-1}-q^{d-2i+2}I)(B-bq^{2i-d}I)
\label{eq:bakb}
\end{eqnarray}
vanishes on $U_i$. Using 
(\ref{eq:keq2}) we find
\begin{eqnarray}
qBK^{-1}-q^{-1}K^{-1}B-b(q-q^{-1})I
\label{eq:qabkb}
\end{eqnarray}
vanishes on $U_i$.
Adding
(\ref{eq:bakb})
to 
$q$ times 
(\ref{eq:qabkb}) 
we find
(\ref{eq:wantkb}) vanishes on $U_i$.
We conclude 
$(K^{-1}-q^{d-2i}I)U_i \subseteq U_{i-1}$ for  $0 \leq i \leq d$.
It follows 
$(K-q^{2i-d}I)U_i \subseteq U_{0}+\cdots + U_{i-1}$ for  $0 \leq i \leq d$.
\\
\noindent
 $\lbrack D^*0\rbrack $:
Let $U_0, U_1, \ldots, U_d$ denote the decomposition
 $\lbrack D^*0\rbrack $.
We show
$KU_i \subseteq U_0+\cdots +U_{i+1}$ for $0 \leq i \leq d$.
Let $i$ be given. We have
\begin{eqnarray*}
KU_i  &\subseteq &
K(U_0+\cdots +U_{i})
\\
&=& K(V^*_{d-i}+\cdots + V^*_d)
\qquad \qquad \qquad (\mbox{by
Lemma 
\ref{thm:decsum}, row $\lbrack D^*0\rbrack $})
\\
&\subseteq& V^*_{d-i-1}+\cdots +V^*_d
\qquad \qquad \qquad (\mbox{by
present Theorem, row
 $\lbrack 0^*D^*\rbrack $})
\\
&=& U_0+\cdots +U_{i+1}
\qquad \qquad  \quad 
\qquad \qquad (\mbox{by
Lemma 
\ref{thm:decsum}, row $\lbrack D^*0\rbrack $}).
\end{eqnarray*}
Next we show
$K^{-1}U_i \subseteq U_{i-1}+\cdots +U_d$ for $0 \leq i \leq d$.
Let $i$ be given. We have
\begin{eqnarray*}
K^{-1}U_i  &\subseteq &
K^{-1}(U_i+\cdots +U_{d})
\\
&=& K^{-1}(V_0+\cdots + V_{d-i})
\qquad \qquad \qquad (\mbox{by
Lemma 
\ref{thm:decsum}, row $\lbrack D^*0\rbrack $})
\\
&\subseteq& V_{0}+\cdots +V_{d-i+1}
\qquad \qquad \qquad (\mbox{by
present Theorem, row
 $\lbrack 0D\rbrack $})
\\
&=& U_{i-1}+\cdots +U_d
\qquad \qquad  \quad 
\qquad \qquad (\mbox{by
Lemma 
\ref{thm:decsum}, row $\lbrack D^*0\rbrack $}).
\end{eqnarray*}
 \\
\noindent
$\lbrack D^*D\rbrack $:
Use the present Theorem, row 
$\lbrack 0^*0\rbrack $ and the involution
given in Remark \ref{rem:inv}(i).
\hfill $\Box $ \\

\noindent We now describe the action of $K^*$ and $K^{*-1}$ on
each of the six decompositions from
Lemma
\ref{thm:sixdecp}.

\begin{theorem}
\label{thm:kkdaction2}
Adopt the assumptions of Definition \ref{def:setup2} and
let $U_0, U_1, \ldots, U_d$ 
denote any one of the six decompositions
of $V$ given in 
Lemma \ref{thm:sixdecp}. Let the map $K^*$
be as in 
Definition \ref{def:b}.
Then for $0 \leq i \leq d$
the action of $K^*$ and $K^{*-1}$ on $U_i$ is described as follows.

\medskip
\centerline{
\begin{tabular}[t]{c|c|c}
      {\rm name} &{\rm action of $K^*$ on $U_i$} & {\rm action of $K^{*-1}$ on $U_i$}
      \\ \hline  \hline
	$\lbrack 0D\rbrack $ &
	$(K^*-q^{d-2i}I) U_i \subseteq U_{i-1}$  &
	$ (K^{*-1}-q^{2i-d}I)U_i\subseteq U_0+\cdots + U_{i-1}$  \\
	$\lbrack 0^*D^*\rbrack $ &
	$ (K^{*}-q^{d-2i}I)U_i\subseteq U_{i+1}+\cdots + U_d$  &
	$(K^{*-1}- q^{2i-d}I)U_i \subseteq U_{i+1}$ \\ 
       $\lbrack 0^*D\rbrack $ & 
	$K^*U_i\subseteq U_{i-1}+\cdots+ U_d$ &
	$K^{*-1}U_i \subseteq U_0+\cdots + U_{i+1}$   \\ 
        $\lbrack 0^*0\rbrack $ & 
	$(K^*-q^{2i-d}I)U_i\subseteq U_{i+1}+\cdots +U_d$ &
	$(K^{*-1}-q^{d-2i}I)U_i \subseteq U_{i+1}$   \\ 
       $\lbrack D^*0\rbrack $ &
        $(K^*-q^{2i-d}I)U_i=0$ &
	$(K^{*-1}-q^{d-2i}I)U_i=0$   \\ 
        $\lbrack D^*D\rbrack $ &
	$(K^*-q^{2i-d}I)U_i\subseteq U_{i-1}$ &
	$(K^{*-1}-q^{d-2i}I)U_i \subseteq U_0+\cdots+U_{i-1}$  
	\end{tabular}}
\medskip
\end{theorem}

\noindent {\it Proof:}
Use Theorem
\ref{thm:kkdaction}
and the involution given in
Remark \ref{rem:inv}(ii).
\hfill $\Box $ \\

\section{The $q$-Serre relations}

\noindent In this section we give two relations involving
the tridiagonal pair $A,A^*$ from Definition
\ref{def:setup2}, and two relations involving
the tridiagonal pair $B,B^*$ from Definition \ref{def:b}.

\begin{theorem}
\label{thm:qsab}
With reference to Definition \ref{def:setup2}
and Definition \ref{def:b},
\begin{eqnarray}
A^3A^*-\lbrack 3 \rbrack_q A^2A^*A+\lbrack 3 \rbrack_q AA^*A^2-A^*A^3&=&0,
\label{eq:aserre}
\\
A^{*3}A-\lbrack 3 \rbrack_q A^{*2}AA^*+\lbrack 3 \rbrack_q A^*AA^{*2}-AA^{*3}&=&0,
\label{eq:asserre}
\\
B^3B^*-\lbrack 3 \rbrack_q B^2B^*B+\lbrack 3 \rbrack_q BB^*B^2-B^*B^3&=&0,
\label{eq:bserre}
\\
B^{*3}B-\lbrack 3 \rbrack_q B^{*2}BB^*+\lbrack 3 \rbrack_q B^*BB^{*2}-BB^{*3}&=&0.
\label{eq:bsserre}
\end{eqnarray}
\end{theorem}
\noindent {\it Proof:}
We first show
(\ref{eq:aserre}).
let $U_0, U_1, \ldots, U_d$ denote the decomposition 
$\lbrack 0D\rbrack $ from Lemma
\ref{thm:sixdecp}.
By Lemma
\ref{thm:aaction}, row $\lbrack 0D\rbrack $,
and using (\ref{eq:eig}), we find that
for $0 \leq i \leq d$
the space $U_i$ is an eigenspace for $A$
with eigenvalue $aq^{2i-d}$.
Abbreviate
$\Psi = 
A^3A^*-\lbrack 3\rbrack_q A^2A^*A+\lbrack 3\rbrack_q AA^*A^2-A^*A^3$.
We show $\Psi =0$.
To do this
we show $\Psi U_i=0$ for $0 \leq i \leq d$.
Let $i$ be given and pick $v \in U_i$.
Observe $A^*v \in 
U_{i-1}+
U_{i}+
U_{i+1}$ by
Lemma
\ref{thm:aaction}, row $\lbrack 0D\rbrack $.
Observe 
$(A-aq^{2i-2-d}I)U_{i-1}=0$,
$(A-aq^{2i-d}I)U_{i}=0$, and
$(A-aq^{2i+2-d}I)U_{i+1}=0$.
By these comments
\begin{eqnarray*}
(A-aq^{2i-2-d}I)
(A-aq^{2i-d}I)
(A-aq^{2i+2-d}I)
 A^*v=0.
\end{eqnarray*}
We may now argue
\begin{eqnarray*}
\Psi v &=&
(A^3A^*-\lbrack 3\rbrack_q A^2A^*A+\lbrack 3\rbrack_q AA^*A^2-A^*A^3)v
\\
&=&
(A^3A^*-\lbrack 3\rbrack_q A^2A^*aq^{2i-d}+
\lbrack 3\rbrack_q AA^*a^2q^{4i-2d}-A^*a^3q^{6i-3d})v
\\
&=& 
(A-aq^{2i-2-d}I)
(A-aq^{2i-d}I)
(A-aq^{2i+2-d}I)A^*v
\\
&=& 0.
\end{eqnarray*}
We have now shown $\Psi U_i=0$ for $0 \leq i \leq d$.
We conclude
 $\Psi =0$ and 
(\ref{eq:aserre}) follows.
To get 
(\ref{eq:asserre}) use
(\ref{eq:aserre}) and the involution 
in Remark \ref{rem:inv}(i).
To get 
(\ref{eq:bserre}),
(\ref{eq:bsserre}) 
apply 
(\ref{eq:aserre}),
(\ref{eq:asserre}) to the tridiagonal pair
$B,B^*$.
\hfill $\Box $ \\


\section{Two modules for  
$U_q(\widehat{ sl}_2)$}

\noindent
In this section we prove the existence part of
Theorem
\ref{thm:main}. We begin with two theorems.

\begin{theorem}
\label{thm:2modaction}
Adopt the assumptions of Definition \ref{def:setup2}.
Let $B,B^*,K$ be as in
Definition \ref{def:k}.
Then
$V$ is an irreducible
$U_q(\widehat{sl}_2)$-module
on which the alternate generators 
act as follows.

\medskip
\centerline{
\begin{tabular}[t]{c|cccccccc}
        {\rm generator}  
       	& $y_0^+$  
         & $y_1^+$  
         & $y_0^-$  
         & $y_1^-$  
         & $k_0$  
         & $k_1$  
         & $k_0^{-1}$  
         & $k_1^{-1}$  
	\\
	\hline 
{\rm action on $V$} 
&  $b^{*-1}B^*$ & $b^{-1}B$ & $a^{*-1}A^*$ & $a^{-1}A$ & $K$ & $K^{-1}$ &
$K^{-1}$ 
& $K$ 
\end{tabular}}
\end{theorem}
\noindent {\it Proof:}
To see that the above action 
on $V$ gives a
$U_q(\widehat{sl}_2)$-module,
compare the equations 
in Theorem 
\ref{thm:key},
Theorem
\ref{thm:kkey},
and Theorem 
\ref{thm:qsab}
with the defining relations for 
$U_q(\widehat{sl}_2)$ 
given in
Theorem \ref{thm:uq2}.
The
$U_q(\widehat{sl}_2)$-module
$V$ is irreducible 
by Definition
\ref{def:tdp}(iv).
\hfill $\Box $ \\

\begin{theorem}
\label{thm:modaction}
Adopt the assumptions of Definition \ref{def:setup2}.
Let $B,B^*,K^*$ be as in
Definition \ref{def:k}.
Then
$V$ is an irreducible
$U_q(\widehat{sl}_2)$-module
on which the alternate generators 
act as follows.

\medskip
\centerline{
\begin{tabular}[t]{c|cccccccc}
        {\rm generator}  
       	& $y_0^+$  
         & $y_1^+$  
         & $y_0^-$  
         & $y_1^-$  
         & $k_0$  
         & $k_1$  
         & $k_0^{-1}$  
         & $k_1^{-1}$  
	\\
	\hline 
{\rm action on $V$} 
&  $a^{-1}A$ & $a^{*-1}A^*$ & $b^{*-1}B^*$ & $b^{-1}B$ & $K^*$ & $K^{*-1}$ &
$K^{*-1}$ & $K^*$ 
\end{tabular}}
\end{theorem}
\noindent {\it Proof:}
To see that the above action 
on $V$ gives a
$U_q(\widehat{sl}_2)$-module,
compare the equations 
in Theorem 
\ref{thm:key},
Theorem
\ref{thm:kkeyd},
and Theorem 
\ref{thm:qsab}
with the defining relations for 
$U_q(\widehat{sl}_2)$ 
given in
Theorem \ref{thm:uq2}.
The
$U_q(\widehat{sl}_2)$-module
$V$ is irreducible 
by Definition
\ref{def:tdp}(iv).
\hfill $\Box $ \\

\noindent It is now a simple matter to prove the existence
part of 
Theorem
\ref{thm:main}.

\medskip
\noindent {\it Proof of Theorem \ref{thm:main} (existence)}:
By Theorem
\ref{thm:2modaction} there exists an irreducible
$U_q(\widehat{sl}_2)$-module structure on $V$
such that $ay_1^-$ acts as $A$ and $a^*y_0^-$
 acts
$A^*$.
By 
Theorem
\ref{thm:modaction} there exists an irreducible
$U_q(\widehat{sl}_2)$-module structure on $V$
such that $ay_0^+$ acts as $A$ and $a^*y_1^+$
 acts as
$A^*$.
\hfill $\Box $ \\

%

\section{Uniqueness} 

\noindent In this section we prove the uniqueness
part of Theorem \ref{thm:main}.

\medskip
\noindent We begin with a comment concerning
finite dimensional irreducible  
$U_q(\widehat{sl}_2)$-modules. 

\begin{lemma}
\label{lem:weight}
Let $V$ denote a finite dimensional irreducible 
$U_q(\widehat{sl}_2)$-module. 
Then there exist
nonzero scalars
$\varepsilon_0,
\varepsilon_1$ in 
$\K$
and 
there exists a decomposition
$U_0, U_1, \ldots, U_d$ of $V$ 
such that both
\begin{eqnarray}
(k_0-\varepsilon_0q^{2i-d}I)U_i=0, \qquad
 (k_1-\varepsilon_1q^{d-2i}I)U_i=0
\qquad 
\qquad (0 \leq i \leq d).
\label{eq:kmove}
\end{eqnarray}
The sequence 
$\varepsilon_0, \varepsilon_1;
U_0, U_1, \ldots, U_d
$
is unique. Moreover for $0 \leq i \leq d$ we have 
\begin{eqnarray}
&&
(\varepsilon_0 y^+_0-q^{d-2i}I)U_i \subseteq U_{i+1}, \qquad
(\varepsilon_1 y^-_1-q^{2i-d}I)U_i \subseteq U_{i+1}, 
\label{eq:altemove1}
\\
&&
(\varepsilon_0 y^-_0-q^{d-2i}I)U_i \subseteq U_{i-1}, \qquad
(\varepsilon_1 y^+_1-q^{2i-d}I)U_i \subseteq U_{i-1}.
\label{eq:altemove2}
\end{eqnarray}
\end{lemma}
\noindent {\it Proof:}
By the construction
$V$ has finite positive dimension.
Since $k_0k_1$
is 
central
in 
$U_q(\widehat{sl}_2)$
and since
$\K$ is algebraically closed,
there exists $\alpha \in \K$ such that
$(k_0k_1-\alpha I)V=0$.
Observe $\alpha \not=0$ since each of $k_0, k_1$ is invertible
on $V$.
For $\theta \in \K$ we define
$V(\theta)=\lbrace v \in V | k_0v=\theta v\rbrace $.
We observe $V(\theta)\not=0$ if and only if
$\theta$ is an eigenvalue of $k_0$ on $V$,
and in this case 
 $V(\theta)$ is the corresponding eigenspace.
For nonzero $\theta \in \K$ we find
using
(\ref{eq:2buq3}), (\ref{eq:2buq4}) that
\begin{eqnarray}
&&
(y_0^+-\theta^{-1}I)
V(\theta) \subseteq V(q^2 \theta), \qquad \qquad
(y_1^--\theta \alpha^{-1}I)
V(\theta) \subseteq V(q^2 \theta),
\label{eq:emovepre1}
\\
&&
(y_0^--\theta^{-1}I)
V(\theta) \subseteq V(q^{-2} \theta), \qquad \qquad
(y_1^+-\theta \alpha^{-1}I)
V(\theta) \subseteq V(q^{-2} \theta).
\label{eq:emovepre2}
\end{eqnarray}
Since $\K$ 
is algebraically closed 
and since $V$ has finite positive dimension, there exists
$\theta \in \K$ such that
$V(\theta)\not=0$.
We observe $\theta \not=0$ since $k_0$ is invertible on $V$.
Since $q$ is not a root of unity the scalars
$\theta, q^{-2}\theta, q^{-4}\theta, \ldots$ are mutually distinct.
These scalars cannot all be eigenvalues of $k_0$ on $V$;
consequently there exists 
a nonzero $\eta \in \K$
such that
$V(\eta)\not=0$ and
$V(q^{-2}\eta)=0$.
Similarly the scalars $\eta, q^2\eta, q^4\eta, \ldots $ are mutually distinct
so they are not all eigenvalues of $k_0$ on $V$;
consequently there exists a nonnegative integer $d$ such
that $V(q^{2i}\eta)$ is nonzero for 
$0\leq i \leq d$ and zero for $i=d+1$. We abbreviate
$U_i=V(q^{2i}\eta)$ for $0 \leq i \leq d$.
From the construction
\begin{eqnarray}
(k_0-q^{2i}\eta I)U_i=0,
\qquad (k_1-\alpha q^{-2i} \eta^{-1}I)U_i=0
\qquad \qquad 
(0 \leq i \leq d).
\label{eq:kk}
\end{eqnarray}
Define 
$\varepsilon_0,
\varepsilon_1$ so that
$\eta=\varepsilon_0 q^{-d}$ and
$\varepsilon_0
\varepsilon_1
=\alpha$.
Observe $\varepsilon_0, \varepsilon_1$ are nonzero.
Eliminating $\eta, \alpha$  in
(\ref{eq:kk}) using
the preceeding equations
we obtain
(\ref{eq:kmove}). 
From (\ref{eq:emovepre1}),
 (\ref{eq:emovepre2}) and our above comments
 we obtain
(\ref{eq:altemove1}),
(\ref{eq:altemove2}), where
$U_{-1}=0$ and $U_{d+1}=0$.
We claim $V=\sum_{i=0}^d U_i$. 
From
(\ref{eq:kmove})--(\ref{eq:altemove2})
we find $\sum_{i=0}^d U_i$ is invariant under each
 of the alternate  generators
for 
$U_q(\widehat{sl}_2)$.
Also $\sum_{i=0}^d U_i$ is nonzero
since each of $U_0, \ldots, U_d$ is nonzero.
We conclude $V=\sum_{i=0}^d U_i$ 
since $V$ is irreducible as a
$U_q(\widehat{sl}_2)$-module.
The sum
$\sum_{i=0}^d U_i$  is direct since each of
$U_0, \ldots, U_d$
is an  eigenspace for $k_0$ and the corresponding eigenvalues
are mutually distinct.
We now see $U_0, \ldots, U_d$ is a decomposition of $V$.
It is clear that
the sequence
$\varepsilon_0, \varepsilon_1;
U_0, U_1, \ldots, U_d
$
is unique.
\hfill $\Box $ \\

\begin{remark}
\rm
We will not use this fact, but
it turns out that the scalars $\varepsilon_0, \varepsilon_1$
from
 Lemma
\ref{lem:weight} are both in $\lbrace 1,-1\rbrace$.
See for example \cite[Proposition 3.2]{cp}. That
proof assumes $\K=\C$ but the assumption is unnecessary.
\end{remark}

\begin{definition}
\label{def:type}
\rm
Referring to
 Lemma
\ref{lem:weight},
we call the sequence 
$U_0, U_1, \ldots, U_d$
the
{\it weight space decomposition} of
 $V$. We call the ordered pair 
 $(\varepsilon_0, \varepsilon_1)$ the {\it type}
 of $V$.
\end{definition}

\begin{example}
Adopt the assumptions of Definition \ref{def:setup2}.
For the $U_q(\widehat{sl}_2)$-module structure  on $V$ 
given in Theorem
\ref{thm:2modaction}
(resp. Theorem 
\ref{thm:modaction}), 
the weight space decomposition 
coincides with
the
decomposition $\lbrack 0^*D \rbrack $ 
(resp. 
$\lbrack D^*0 \rbrack $)
from
Lemma
\ref{thm:sixdecp}.
Both 
module structures  have type $(1,1)$.
\end{example}
\noindent {\it Proof:}
We first consider the
$U_q(\widehat{sl}_2)$-module structure
from
 Theorem
\ref{thm:2modaction}.
Let $U_0, U_1, \ldots, U_d$ denote the
decomposition $\lbrack 0^*D \rbrack $.
 By Definition
\ref{def:k}(iii) we find $(K-q^{2i-d}I)U_i=0$ for
$0 \leq i\leq d$. 
By Theorem
\ref{thm:2modaction} we find
$k_0, k_1$ act on $V$ as 
 $K, K^{-1}$ respectively.
 Therefore
$(k_0-q^{2i-d}I)U_i=0$ 
and 
$(k_1-q^{d-2i}I)U_i=0$ 
for
$0 \leq i \leq d$.
Define $\varepsilon_0=1$,
$\varepsilon_1=1$ and observe these values
satisfy
(\ref{eq:kmove}).
By Definition \ref{def:type}, $V$
has weight space decompostion
$U_0, U_1, \ldots, U_d$ and
 type $(1,1)$.
We have now proved our assertions concerning
the $U_q(\widehat{sl}_2)$-module structure
from
 Theorem
\ref{thm:2modaction}.
The proof for the
$U_q(\widehat{sl}_2)$-module structure
from
 Theorem
\ref{thm:modaction}  is similar and omitted.
\hfill $\Box $ \\

\medskip
\noindent {\it Proof of Theorem \ref{thm:main}(uniqueness)}:
For $0 \leq i \leq d$ let $V_i$ (resp. $V^*_i$) denote
the eigenspace of $A$ (resp. $A^*$) associated with
$\theta_i$ (resp. $\theta^*_i$).
We assume 
a $U_q(\widehat{sl}_2)$-module structure on
$V$ such that
$ay_1^{-}$ acts as $A$
and 
$a^*y_0^{-}$ acts as $A^*$.
We show
the alternate generators for
$U_q(\widehat{sl}_2)$
act on $V$
according to the table of Theorem 
\ref{thm:2modaction}. Observe
the 
$U_q(\widehat{sl}_2)$-module structure is irreducible
in view of Definition
\ref{def:tdp}(iv).
Let 
$(\varepsilon_0, \varepsilon_1)$ 
denote the type
of the 
$U_q(\widehat{sl}_2)$-module structure.
We claim 
$(\varepsilon_0, \varepsilon_1)=
(1,1)$.
To see this,
consider the weight space decomposition
$U_0, U_1, \ldots, U_d$
from Lemma \ref{lem:weight}.
By 
(\ref{eq:altemove1})
and since $ay_1^-$ acts on $V$ as $A$
we find
\begin{eqnarray}
\label{eq:aactionfirst}
(\varepsilon_1A- a q^{2i-d}I)U_i \subseteq U_{i+1}
\qquad \qquad (0 \leq i \leq d).
\end{eqnarray}
Similarly 
\begin{eqnarray}
\label{eq:asactionfirst}
(\varepsilon_0 A^*- a^* q^{d-2i}I)U_i \subseteq U_{i-1}
\qquad \qquad (0 \leq i \leq d).
\end{eqnarray}
From 
(\ref{eq:aactionfirst}) 
 we find that
for $0 \leq i \leq d$ the scalar $\varepsilon^{-1}_1 a q^{2i-d}$
is an eigenvalue of $A$  and the dimension
of the corresponding eigenspace has the same dimension as
$U_i$. Apparently the sequence
$\varepsilon^{-1}_1 a q^{2i-d}$ $(0 \leq i \leq d)$
is an ordering the eigenvalues of $A$.
Recall
$\theta_i= a q^{2i-d}$
for $0 \leq i \leq d$.
Therefore the sequence 
$\varepsilon^{-1}_1 a q^{2i-d}$ $(0 \leq i \leq d)$
is a permutation of the sequence
$ a q^{2i-d}$ $(0 \leq i \leq d)$.
Since $q$ is not a root of unity we must have
$\varepsilon_1=1$.
By a similar argument we find
$\varepsilon_0=1$.
Setting
$(\varepsilon_0, \varepsilon_1)=(1,1)$ in
(\ref{eq:aactionfirst}), 
(\ref{eq:asactionfirst}) 
we find 
$(A-\theta_iI)U_i \subseteq U_{i+1} $
and 
$(A^*-\theta^*_iI)U_i \subseteq  U_{i-1}$
for 
$0 \leq i \leq d$.
By this and
 \cite[Theorem  4.6]{TD00}
we find
$U_i = (V^*_0+\cdots + V^*_i)\cap (V_i + \cdots + V_d)$
for $0 \leq i\leq d$.
In other words $U_0,  \ldots, U_d$
is the decomposition $\lbrack 0^*D\rbrack $ from
Lemma
\ref{thm:sixdecp}.
By Definition \ref{def:k}(iii) we have 
$(K-q^{2i-d}I)U_i=0$ for $0 \leq i \leq d$.
Comparing this with 
(\ref{eq:kmove}) and recalling 
$(\varepsilon_0, \varepsilon_1)=(1,1)$
we find $k_0, k_1$ act on $V$ as $K, K^{-1}$ respectively.
 Apparently
 $k^{-1}_0, k^{-1}_1$ act on $V$ as $K^{-1}, K$ respectively.
We show $by_1^+$ acts on $V$ as $B$.
Define $W=\lbrace v \in V | (by^+_1-B)v=0 \rbrace $.
We show $W=V$.
To do this we show $W\not=0$, $AW\subseteq W$,
 $A^*W\subseteq W$.
Observe $(B-bq^{-d}I)U_0=0$ by
Theorem \ref{thm:baction}, row $\lbrack 0^*D\rbrack $.
Observe $(y^+_1-q^{-d}I)U_0=0$ by
(\ref{eq:altemove2}). By these comments
$by_1^+-B$ vanishes on $U_0$. Therefore 
$U_0 \subseteq W$ so
$W\not=0$.
By 
(\ref{eq:2buq5}) (with $i=1$),
by (\ref{eq:eq1}), and since $ay^{-}_1, A$ agree on $V$,
we find
$(by^+_1 -B)A$, $q^2 A(by^+_1 -B)$ agree on $V$.
Using this we find
$AW\subseteq W$.
By 
(\ref{eq:2buq6}) (with $i=1$),
by (\ref{eq:eq2}), and since $a^*y^{-}_0, A^*$ agree on $V$,
we find
$(by^+_1 -B)A^*$, $q^{-2} A^*(by^+_1 -B)$ agree on $V$.
Using this we find
$A^*W\subseteq W$.
We have now shown
$W\not=0$, $AW\subseteq W$, $A^*W\subseteq W$. 
Now $W=V$ in view of
Definition
\ref{def:tdp}(iv).
We conclude $(by_1^+-B)V=0$ so $by_1^+$ acts on $V$ as $B$.
By a similar argument we find  
$b^*y_0^+$ acts on $V$ as $B^*$.
We have now shown
$y^{\pm}_i$, $k_i^{{\pm}1}$, $i\in \lbrace 0,1\rbrace $
act on $V$ according to the table of
Theorem \ref{thm:2modaction}. 
It follows the given
$U_q(\widehat{sl}_2)$-module structure 
is unique.
By a similar argument we obtain the uniqueness of
 the irreducible $U_q(\widehat{sl}_2)$-module structure on $V$
such that
$ay^{+}_0$ acts as $A$ and 
$a^*y^{+}_1$ acts as $A^*$.
\hfill $\Box $ \\

\section{Comments}

\noindent We have a comment on
Theorem
\ref{thm:main}.

\begin{lemma}
\label{lem:whenir}
Let $a, a^*$ denote nonzero scalars in $\K$.
Let $A,A^*$ denote elements in
$U_q(\widehat{sl}_2)$ which satisfy
\begin{eqnarray}
A = ay_1^-, \qquad \qquad 
A^* = a^*y_0^- 
\label{eq:aas}
\end{eqnarray}
or 
\begin{eqnarray}
A = ay_0^+, \qquad \qquad 
A^* = a^*y_1^+.
\label{eq:2aas}
\end{eqnarray}
Let $V$ denote a finite dimensional irreducible
$U_q(\widehat{sl}_2)$-module of type $(1,1)$.
Assume $V$ is irreducible as an $(A,A^*)$-module.
Then the pair $A,A^*$ acts on $V$ as a tridiagonal pair.
Denoting the diameter of this pair by $d$,
the sequence $aq^{2i-d}$ $(0 \leq i \leq d)$ is
a standard ordering of the eigenvalues for $A$
on $V$ and 
the sequence $a^*q^{d-2i}$ $(0 \leq i \leq d)$ is
a standard ordering of the eigenvalues for $A^*$ on $V$.
\end{lemma}
\noindent {\it Proof:}
First assume (\ref{eq:aas}). By 
(\ref{eq:2buq7}) 
and
(\ref{eq:aas}) 
we find
both
\begin{eqnarray}
A^3A^*-\lbrack 3 \rbrack_q A^2A^*A+\lbrack 3 \rbrack_q AA^*A^2-A^*A^3&=&0,
\label{eq:a3}
\\
A^{*3}A-\lbrack 3 \rbrack_q A^{*2}AA^*+
\lbrack 3 \rbrack_q A^*AA^{*2}-AA^{*3}&=&0.
\label{eq:as3}
\end{eqnarray}
Let $U_0, U_1, \ldots, U_d$ denote the
weight space decomposition of $V$ from
Definition
\ref{def:type}.
Setting $(\varepsilon_0, \varepsilon_1)=(1,1)$
in
Lemma
\ref{lem:weight}
and using
(\ref{eq:aas})
we find
both
\begin{eqnarray}
(A-aq^{2i-d}I)U_i &\subseteq& U_{i+1}  \qquad \qquad (0 \leq i \leq d),
\label{eq:incl1}
\\
(A^*-a^*q^{d-2i}I)U_i &\subseteq& U_{i-1}
\qquad \qquad (0 \leq i \leq d).
\label{eq:incl2}
\end{eqnarray}
We draw several conclusions from these lines.
From 
(\ref{eq:incl1})
(resp. 
(\ref{eq:incl2}))
the action of $A$ (resp. $A^*$) 
on $V$
is diagonalizable.
Also for $0 \leq i \leq d$ the scalar $aq^{2i-d}$
(resp. $a^*q^{d-2i}$)
is an eigenvalue for this action and the 
corresponding 
eigenspace has the same dimension as
$U_i$.
In particular the scalars $aq^{2i-d}$ $(0 \leq i \leq d)$ 
(resp.  $a^*q^{d-2i}$ $(0 \leq i \leq d)$)
are the eigenvalues of $A$ (resp. $A^*$) on $V$.
We are assuming  $V$ is irreducible as an $(A,A^*$)-module.
This means there does not exist a 
subspace $W\subseteq V$ such
that
$AW\subseteq W$,
$A^*W\subseteq W$,
$W\not=0$,  
$W\not=V$.
We show $A,A^*$
acts on $V$ as a tridiagonal pair. To do this
we apply \cite[Example 1.7]{TD00}.
In order to apply this example we must show neither of
$A,A^*$ is nilpotent on $V$.
We mentioned above that
each of $A,A^*$ is diagonalizable on $V$. Neither of
$A,A^*$ is zero on $V$ so neither of $A,A^*$ is nilpotent
on $V$.
Now by \cite[Example 1.7]{TD00}
we find $A,A^*$ act on
$V$ as a tridiagonal pair.
The diameter of this pair is $d$ since
each of $A, A^*$ has $d+1$ distinct eigenvalues.
By \cite[Lemma 4.8]{qSerre}
 there exists a standard ordering of the eigenvalues
of $A$ (resp. $A^*$) on $V$ of the form
$\alpha q^{2i-d}$ $(0 \leq i \leq d)$
(resp. $\alpha^*q^{d-2i}$ $(0 \leq i \leq d))$,
where $\alpha$ (resp. $\alpha^*$) is an
appropriate nonzero scalar in $\K$.
Combining this with our above remarks
we find 
$\alpha=a$ and $\alpha^*=a^*$.
Therefore the sequence
$a q^{2i-d}$ $(0 \leq i \leq d)$ is a standard ordering
of the eigenvalues for $A$ on $V$
and the sequence
$a^*q^{d-2i}$ $(0 \leq i \leq d)$ is
a standard ordering of the eigenvalues for $A^*$ on $V$.
We have now proved the result for case
(\ref{eq:aas}). For the case   
(\ref{eq:2aas}) the proof is similar and omitted.
\hfill $\Box $ \\


\section{Suggestions for further research}

\noindent In this section we give some open problems. 
The first problem is motivated
by Lemma
\ref{lem:whenir}.

\begin{problem}
\label{prob:one}
\rm
Let $a,a^*$ denote nonzero scalars in $\K$ and 
let $A,A^*$ denote the elements of
$U_q(\widehat{sl}_2)$ given in
(\ref{eq:aas}) or
(\ref{eq:2aas}).
 Let $V$ denote a finite dimensional
irreducible 
$U_q(\widehat{sl}_2)$-module of type $(1,1)$.
Find a necessary and sufficient condition for
$V$ to be irreducible as an $(A,A^*)$-module.
\end{problem}

\noindent In order to state the next problem
we recall a few terms. Let
$V$ denote a vector space over $\K$ with finite positive
dimension. Let $\mbox{End}(V)$ denote
the $\K$-algebra consisting of all linear transformations
from $V$ to $V$. By an 
{\it antiautomorphism} of
 $\mbox{End}(V)$ we mean a 
$\K$-linear bijection $\dagger : 
 \mbox{End}(V) \rightarrow
 \mbox{End}(V)$
such that $(XY)^\dagger = Y^\dagger X^\dagger $
for all $X, Y \in 
 \mbox{End}(V)$.

\begin{problem} 
\rm
Let $A,A^*$ denote a tridiagonal pair
on $V$. Show there exists an antiautomorphism $\dagger$ of 
 $\mbox{End}(V)$ such that $A^\dagger = A$ and
$A^{*\dagger}=A^*$.
We remark that $\dagger $ exists if $A,A^*$ is a Leonard pair
\cite[Theorem 7.1]{qrac}.
\end{problem}

\section{Acknowledgements}
The second author
would like to thank Georgia Benkart for pointing out
around 1997 that the two mysterious equations
that were showing up in connection with tridiagonal pairs
are known to researchers in quantum groups
as the $q$-Serre relations.
Both authors would like to thank 
Kenichiro Tanabe for giving us a week-long
tutorial in the summer of 1999 on
the subject of
$U_q(\widehat{sl}_2)$ and its modules;
the resulting boost in our
understanding illuminated the way to this paper.

\noindent Tatsuro Ito \hfil\break
\noindent Department of Computational Science \hfil\break
\noindent Faculty of Science \hfil\break
\noindent Kanazawa University \hfil\break
\noindent Kakuma-machi \hfil\break
\noindent Kanazawa 920-1192, Japan \hfil\break
\noindent email:  {\tt ito@kappa.s.kanazawa-u.ac.jp}

\bigskip

\noindent Paul Terwilliger \hfil\break
\noindent Department of Mathematics \hfil\break
\noindent University of Wisconsin \hfil\break
\noindent Van Vleck Hall \hfil\break
\noindent 480 Lincoln Drive \hfil\break
\noindent Madison, WI 53706-1388 USA \hfil\break
\noindent email: {\tt terwilli@math.wisc.edu }\hfil\break

\end{document}